%

\input amstex
\documentstyle{amsppt}
\magnification=1200
 \vsize19.5cm
  \hsize13.5cm

   \TagsOnRight

\pageno=1
\baselineskip=15.0pt
\parskip=3pt



\def\p{\partial}
\def\noo{\noindent}
\def\eps{\varepsilon}

\def\Om{\Omega}

\def\pom{{\p \Om}}
\def\bom{{\overline\Om}}
\def\R{\bold R}

\def\th{\theta}
\def\wtt{\widetilde}

\def\Ga{\Gamma}

\def\diam{\text{diam}}
\def\dist{\text{dist}}
\def\det{\text{det}}

\def\ol{\overline}
\def\lan{\langle}
\def\ran{\rangle}
\def\D{\nabla}
\def\phi{\varphi}
\def\C{\Cal C}

\def\M{\Cal M}
\def\N{\Cal N}

\nologo
\NoRunningHeads

\topmatter

\title  {The affine Plateau problem}\endtitle

\author{\bf Neil S. Trudinger \ \ \ \ Xu-Jia Wang}\endauthor

\affil{Centre for Mathematics and Its Applications\\
       The Australian National University\\
       Canberra, ACT 0200\\
       Australia}\endaffil

\address{Centre for Mathematics and Its Applications,
       Australian National University, \newline
       Canberra, ACT 0200,
       Australia}\endaddress

\email{Neil.Trudinger\@maths.anu.edu.au,\ \ \ \ \ \
X.J.Wang\@maths.anu.edu.au}\endemail

\thanks{Supported by Australian Research Council}\endthanks

\abstract {In this paper, we study a second order variational
problem for locally convex hypersurfaces, which is the affine
invariant analogue of the classical Plateau problem for minimal
surfaces. We prove existence, regularity and uniqueness results
for hypersurfaces maximizing affine area under appropriate
boundary conditions.}\endabstract

\endtopmatter

\document

\baselineskip=14.0pt
\parskip=2.5pt

\centerline {\bf \S 1. Introduction}

\vskip10pt

In this paper we study the Plateau problem for affine maximal
hypersurfaces, which is the affine invariant analogue of the
classical Plateau problem for minimal surfaces. In particular we
formulate the affine Plateau problem as a geometric variational
problem for the affine area functional, and prove the existence
and regularity of maximizers. As a special case, we obtain
corresponding existence and regularity results for the variational
Dirichlet problem for the fourth order affine maximal surface
equation, together with a uniqueness result for generalized
solutions.

The affine Plateau problem may be formulated as follows. Let
$\M_0$ be a bounded, connected hypersurface in Euclidean
$(n+1)$-space, $\R^{n+1}$, with smooth boundary $\Ga=\p\M_0$.
Assume that $\M_0\cup \Ga$ is smooth and locally uniformly convex
up to boundary. Let $S[\M_0]$ denote the set of locally uniformly
convex hypersurfaces $\M$ with boundary $\Ga$, which can be
smoothly deformed from $\M_0$ in the family of locally uniformly
convex hypersurfaces whose Gauss mapping images lie in that of
$\M_0$. A hypersurface $\M\subset\R^{n+1}$ is called {\it locally
uniformly convex} if it is a $C^2$ immersion of an $n$-manifold
$\N$, whose principal curvatures are everywhere positive, and the
{\it Gauss mapping} of $\M$ is the mapping $G:\ \M\to S^n$ which
assigns to every point in $\M$ its unit normal vector. The {\it
affine metric} (also called the Berwald-Blaschke metric) on a
locally uniformly convex hypersurface $\M$ is defined by
$$g=K^{-1/(n+2)}I\!I , \tag 1.1$$
where $K$ is the Gauss curvature and $I\!I$ is the second
fundamental form of $\M$. From the affine metric we have the
affine area functional
$$A(\M)=\int_\M K^{1/(n+2)}. \tag 1.2$$
The affine Plateau problem is to determine a hypersurface $\M\in
S[\M_0]$, maximizing the functional $A$ over $S[\M_0]$, that is
$$A(\M)=\sup\{ A(\M')\ {|}\ \M'\in S[\M_0]\}.\tag 1.3$$
Recall that the classical Plateau problem is formulated as a
variational minimization problem for the Euclidean area. The two
dimensional affine Plateau problem was raised by Chern in his
pioneering article [9], see also Calabi [7]. For this we have the
following solution; (Theorems 3.1, 6.2, and 8.3).

\proclaim{Theorem A} There exists a smooth, locally uniformly
convex hypersurface $\M\in S[\M_0]$ solving the variational
Plateau problem (1.3) in the two dimensional case, $n=2$, if and
only if the image of the Gauss mapping of $\M_0$ does not cover
any hemisphere.
\endproclaim

In higher dimensions we shall prove that the affine Plateau
problem is solvable in a generalized sense (Theorem 3.1), with
solutions locally uniformly convex and smooth if they are locally
strictly convex (Theorem 6.2).

A special case of the affine Plateau problem occurs when $\M_0$ is
the graph of a smooth, locally uniformly convex function $\phi$
defined on the closure of a bounded smooth domain
$\Om\subset\R^n$. If $\M=\M_u$ is the graph of a convex function
$u\in C^2(\Om)$, the affine area of $\M$ is given by
$$A(u):=A(\M_u)=\int_\Om [\det D^2 u]^{1/(n+2)}. \tag 1.4$$
The functional $A$ in (1.4) is concave and upper semi-continuous
with respect to local uniform convergence. If a smooth, locally
uniformly convex function $u$ is stationary for $A$, then $u$
satisfies the {\it affine maximal surface equation},
$$ L[u] := U^{ij}w_{ij}=0,\tag 1.5$$
where
$$w =[\det D^2 u]^{-(n+1)/(n+2)} ,\tag 1.6$$
and $[U^{ij}]$ is the cofactor matrix of the Hessian matrix $D^2
u$, which is divergence free for any fixed $i$ or $j$. The
subscripts $i, j$ denote partial derivatives with respect to the
variables $x_i, x_j$. $L$ is a nonlinear fourth order partial
differential operator, which is elliptic with respect to locally
uniformly convex solutions. The concavity of $A$ implies that the
affine maximal surface equation (1.5) is both a necessary and
sufficient condition for a smooth, locally uniformly convex
function $u$ to locally maximize the affine area functional (1.4).

In the graph case the set $S[\M_0]=S[\phi]\  (=S[\phi, \Om])$
consists of locally uniformly convex functions $u\in C^2(\Om)\cap
C^0(\bom)$ satisfying $u=\phi$ on $\pom$ and $Du(\Om)\subset
D\phi(\Om)$. Accordingly we have the variational problem of
finding a function $u\in S[\phi]$, maximizing the affine area
functional $A$ on $S[\phi]$, that is
$$A(u)=\sup \{A(v)\ {|} \ v\in S[\phi]\} ,\tag 1.7$$
and as a special case of Theorem A, this problem is solvable for
$n=2$. More generally, we will extend the definition of affine
area to non-smooth locally convex functions (and locally convex
hypersurfaces) in Section 2 and study the variational problem of
maximizing the extended affine area $A$ on $\ol S[\phi]$, the
closure of $S[\phi]$ under uniform convergence, that is to find
$u\in\ol S[\phi]$ such that
$$A(u)=\sup \{A(v)\ {|} \ v\in \ol S[\phi]\} .\tag 1.8$$
We prove the existence of a {\it unique} maximizer $u$ for the
extended functional (1.8), which is smooth if it is locally
strictly convex; (Theorems 2.1 and 6.2)

\proclaim{Theorem B} There exists a unique locally convex function
$u$ solving the variational boundary value problem (1.8), in all
dimensions, which is smooth and locally uniformly convex in the
interior of any set where it is locally strictly convex.
\endproclaim

The variational problems (1.7) and (1.8) extend the {\it first
boundary value problem} for the affine maximal surface equation
(1.5),
$$\eqalign{
 u & = \phi \ \ \text{on}\ \ \pom,   \cr
 Du & = D\phi\ \ \text{on}\ \ \pom,  \cr} \tag 1.9$$
for if we have a classical, locally uniformly convex solution
$u\in C^4(\Om)\cap C^1(\bom)$ of (1.5), (1.9), $u$ will also solve
(1.7) uniquely. Note that if the domain $\Om$ is convex, a locally
convex function is convex, namely its graph lies above its tangent
planes everywhere.

The existence of at least non-smooth maximizers in Theorems A and
B follows from the upper semi-continuity of the affine area
functional (1.4). For Theorem A we also need the fundamental lemma
from [22] that a locally convex hypersurface with boundary on a
hyperplane is convex, which facilitates the use of uniform local
graph representations. The existence and uniqueness of the
maximizers for the variational problem (1.8) will be proved in
Section 2, (Theorem 2.1). For the uniqueness we need a preliminary
result, Lemma 2.3, which guarantees that the Monge-Amp\`ere
measure of a maximizer is regular. The general Plateau problem is
treated in Section 3, where we prove the existence of maximizers
for the corresponding extended variational problem (1.3) in the
class of locally convex hypersurfaces, (Theorem 3.1), under the
necessary condition that the Gauss map of $\M_0$ does not cover
any hemi-sphere.

The rest of this paper, Sections 4 to 8, deals with the issue of
regularity. In Section 4 we reprove and extend to inhomogeneous
equations
$$L[u]=f, \tag 1.10$$
a priori estimates for smooth, locally uniformly convex solutions,
already established in [21]; (Theorems 4.1 and 4.2). These
estimates are then employed to obtain corresponding regularity
results through approximation by solutions to the {\it second
boundary value problem},
$$ \eqalign
{ u & =\phi \ \ \ \text{on}\ \ \pom,\cr
  w & =\psi \ \ \ \text{on}\ \ \pom, \cr  }
          \tag 1.11 $$
which is the Dirichlet problem when (1.5) (1.6) is considered as a
second order system of two equations in $u$ and $w$. In Section 5
we prove the existence of locally smooth solutions to the boundary
value problem (1.10), (1.11), (Theorem 5.1). Using a penalty
method, we then prove in Section 6 that maximizers of the extended
functional (1.8) can be approximated locally in $\Om$ by smooth,
uniformly convex solutions of the affine maximal surface equation
(1.5), (Theorem 6.1). As consequences, we conclude a fundamental
regularity result in all dimensions, that strictly convex affine
maximal functions are smooth, (Theorem 6.2), thereby completing
the proof of Theorem B, and extend our two dimensional Bernstein
Theorem in [21] to general convex affine maximal functions,
(Theorem 6.3).

As with the Monge-Amp\`ere equation, the strict convexity of
solutions is a critical issue for the affine maximal surface
equation. In [21] we proved the Bernstein theorem holds in any
dimension under an asymptotic, strict convexity assumption, which
is automatically satisfied if any affine maximal hypersurface with
boundary on a hyperplane is strictly convex, which we proved for
affine maximal surfaces in $\R^3$. In this paper we prove the
maximizers in Theorems A and B are strictly convex for $n=2$ under
general boundary conditions. The argument in this case becomes
much more complicated than that in [21], as we have to eliminate
the possibility of straight line segments in graphs having both
endpoints at the boundary. In Section 7 we treat the extension of
the Legendre transform to locally convex functions in arbitrary
domains in preparation for our treatment of the strict convexity
in Section 8. There we prove that,  in two dimensions the
maximizers of the affine area functionals (1.3) and (1.8) are
locally   strictly convex,  thereby completing the proof of
Theorem A, (Theorems 8.2 and 8.3). Finally we extend the local
strict convexity result to the inhomogeneous equation (1.10),
(Theorem 8.4), yielding the unconditional interior regularity of
maximizers in the two dimensional case.

The inhomogeneous equation (1.10), which is crucial for our
approximation arguments, is a {\it prescribed affine mean
curvature} equation, as the quantity
$$H_A[u]=-\frac {1}{n+1} L[u]\tag 1.12$$
is the {\it affine mean curvature} of the graph of $u$ [1, 16,
20], and an affine maximal graph is one with vanishing affine mean
curvature. The corresponding variational problem to maximize the
functional
$$A_f(u)=A(u)-\int fu \tag 1.13$$
is actually treated in Section 2 for bounded $f$, with the
inhomogeneous version of Theorem B proved in Theorems 2.1 and 6.2.
From our last regularity result, Theorem 8.4, we obtain the full
interior regularity of maximizers in the two dimensional case.

We conclude the paper with some remarks on the boundary behavior
of maximizers and the regularity problem in higher dimensions.

\vskip30pt

\baselineskip=14.50pt
\parskip=2.5pt

\centerline{\bf \S 2. The affine area functional}

\vskip10pt

In this section,  we prove the existence and uniqueness of
maximizers for the variational problem (1.8). First we need to
extend the affine area functional (1.4) to arbitrary locally
convex functions and prove its upper semi-continuity. Note that as
we will be including domains which are not convex, it is necessary
to distinguish locally convex functions from convex functions.
However, when there is no ambiguity, we will typically abbreviate
locally convex to convex.

We begin with the definition of the Monge-Amp\`ere measure. Let
$\Om$ be a bounded domain in $\R^n$ and $u$ a convex function in
$\Om$. The normal mapping of $u$, $N_u$, is a multi-valued
mapping, defined as follows [17]. For any point $x\in \Om$,
$N_u(x)$ is the set of slopes of support hyperplanes of $u$ at
$x$, that is
$$N_u(x)= \{p\in\R^n\ {|}\ u(y)\ge u(x)+p\cdot (y-x)
     \ \forall\ y\in\Om\};$$
and for any Borel set $E$, $N_u(E)=\bigcup_{x\in E} N_u(x)$. When
$u$ is $C^1$, the normal mapping coincides with the gradient
mapping $Du$.

From the normal mapping one introduces the Monge-Amp\`ere measure
$\mu[u]$, which is a Radon measure given by
$$\mu[u](E)= |N_u(E)|, $$
that is the Lebesgue measure of $N_u(E)$.  When $u$ is $C^2$ we
have $\mu[u]=(\det D^2 u) dx$.

It is a basic result of Aleksandrov that $\mu[u]$ is weakly
continuous with respect to the convergence of convex functions.
That is if $\{u_j\}$ is a sequence of convex functions converging
to $u_0$ in $L^1_{loc}$, then $\mu[u_j]$ converges weakly to
$\mu[u_0]$ as measures. It follows that for any closed set
$E\subset\Om$,
$$\mu[u_0](E)\ge \lim_{j\to\infty} \mu[u_j] (E).\tag 2.1 $$

The Monge-Amp\`ere measure $\mu[u]$, as a Radon measure, can be
decomposed as the sum of a regular part and a singular part. That
is
$$\mu[u]=\mu_r[u] +\mu_s[u],$$
where the regular part $\mu_r[u]$ is a measure given by a locally
integrable function, $\mu_r[u]= (\mu_r[u])dx$ (we use the same
notation $\mu_r[u]$ to denote the Radon-Nikodym derivative of the
measure $\mu_r[u]$) and, the singular part $\mu_s[u]$ is a measure
supported on a set with Lebesgue measure zero.

The regular part $\mu_r[u]$ is determined by the function $u$
explicitly. Indeed, since $u$ is convex, it is twice
differentiable almost everywhere. In this paper we will use the
notation
$$(\p^2 u)=(\p_{ij}u)$$
to denote the Hessian matrix of a function if it is twice
differentiable almost everywhere. Then $\det \p^2 u$ is a
measurable function. We have $\p^2 u=D^2 u$ when $u$ is $C^2$
smooth.

\proclaim {Lemma 2.1} Let $u$ be a convex function. Then
 $$\mu_r[u] = \det\, \p^2 u. \tag 2.2$$
\endproclaim

\noo {\it Proof}.  Let $u_h$ be the mollification of $u$. That is
$$u_h(x)=h^{-n} \int_\Om u(y) \eta(\frac{x-y}{h}) dy \tag 2.3$$
for some nonnegative smooth function $\eta$ supported on the unit
ball $B_1(0)$ and satisfying $\int_{\R^n} \eta=1$, where
$x\in\Om_h=\{x\in\Om\ {|}\ \dist(x, \pom)>h\}$. Then at any point
$x\in\Om$ where $u$ is twice differentiable, we have $D^2 u_h(x)
\to \p^2 u(x)$ [25]. Hence for any measurable set $E\subset\subset
\Om$,
$$\int_E \det \p^2 u
    \le \lim_{h\to 0}\int_E \det D^2 u_h . \tag 2.4$$
By the weak convergence (2.1), we have for any closed set
$E\subset \Om$ such that $\mu_s[u](E)=0$,
$$\int_E \det \p^2 u\le \mu_r[u](E).$$
That is
$$ \det \p^2 u \le \mu_r[u]\ \ \ a.e.. \tag 2.5$$
We claim that the equality holds in (2.5). Indeed, for a.e.
$x_0\in\Om$ we have,
$$\align
\mu_r[u] (x_0)
& =\lim_{\eps\to 0} \frac {\mu_r[u](B_\eps (x_0))}{|B_\eps(x_0)|}\\
& \le \lim_{\eps\to 0} \frac {\mu [u](B_\eps (x_0))}{|B_\eps(x_0)|}.\\
\endalign
$$ If $u$ is twice differentiable at $x_0$, then for any $x\in
B_\eps(x_0)$,
$$|Du(x)-Du(x_0)-\p^2u(x_0) (x-x_0)|\le \delta |x-x_0|$$
for some constant $\delta>0$, with $\delta\to 0$ as $\eps\to 0$.
Suppose for simplicity that $Du(x_0)=0$. Let $w=\frac 12
x\cdot(\p^2 u(x_0)+\delta I) \cdot x$, where $I$ is the unit
matrix. Then for $|x-x_0|$ sufficiently small,
$N_u(B_\eps(x_0))\subset N_w(B_\eps(x_0))$. It follows
$$\align
\mu_r[u] (x_0)
 & \le  \lim_{\eps\to 0}
        \frac {|N_u(B_\eps(x_0))|}{|B_\eps(x_0)|}\\
 & \le  \lim_{\eps\to 0}\frac {|N_w(B_\eps(x_0))|}{|B_\eps(x_0))|}\\
 & =  \det (\p^2 u(x_0)+\delta I)\\
 & \le  \det (\p^2 u(x_0))+C\delta ,\\
\endalign $$
where $C$ depends on $\p^2 u(x_0)$ but is independent of $\eps$.
Sending $\eps$ to zero, we obtain $\mu_r[u] \le\det \p^2 u$ at
$x_0$. Hence we have $\mu_r[u]= \det \p^2 u$. $\square$

Similarly for a given (non-smooth) convex hypersurface $\M$ one
can introduce the (Gauss) curvature measure on $\M$. Let $G: \M\to
S^n$ denote the generalized Gauss mapping. That is for any point
$p\in \M$, $G(p)$ is the set of normals of the support hyperplanes
of $\M$ at $p$. The curvature measure $\mu[\M]$ is defined, for
any Borel set $E\subset\M$, by
$$\mu[\M] (E)= |G(E)| . \tag 2.6$$
The curvature measure can also be decomposed as the sum of a
regular part and a singular part, namely
$$\mu[\M]=\mu_r[\M] +\mu_s[\M] . $$
A similar proof as that of Lemma 2.1 shows that the regular part
$\mu_r[\M]$ is given by the Gauss curvature of $\M$, which is well
defined almost everywhere. In particular if $\M$ is the graph of a
convex function $u$, then
$$\mu_r[\M]=\frac {\det \p^2 u}{(1+|Du|^2)^{(n+2)/2}} . \tag 2.7$$

\vskip10pt

 By Lemma 2.1, we can extend the definition of the affine area
functional (1.4) from smooth convex functions to general convex
functions by
$$A(u) = A(u, \Om) = \int_\Om [\det \p^2 u]^{1/(n+2)}. \tag 2.8$$
By (2.7) we have
$$A(\M)=\int_\M \big(\mu_r [\M]\big)^{1/(n+2)} . \tag 2.9$$
Formulae (2.8) (2.9) imply that the affine area is invariant under
unimodular affine transformations in $\R^{n+1}$, and in particular
is independent of the choice of the coordinates. We remark that
different but equivalent definitions for the affine surface area
have been introduced [14]. Our definition (2.8) is more
straightforward.

\proclaim {Definition} A (locally) convex function $u$ in a domain
$\Om$ is called affine maximal if it is a maximizer of (2.8) under
local convex perturbation. That is for any (locally) convex
function $v$ such that $u-v$ has compact support in
$\Om'\subset\subset\Om$, $A(v, \Om')\le A(u, \Om')$. A locally
convex hypersurface $\M$ is called affine maximal if locally it is
the graph of an affine maximal function.
\endproclaim

Since the function $r\to (\det\, r)^{1/(n+2)}$ is strictly concave
on the cone of positive symmetric $n\times n$ matrices, so also is
the functional $A$, that is
$$ A(tu+(1-t)v)\ge tA(u)+(1-t)A(v)\tag 2.10$$
for all convex $u, v$ and $0\le t\le 1$, with equality holding if
and only if $\p^2 u=\p^2 v$.  Moreover if $u$ is a convex function
and $\eta$ is a continuous function twice differentiable a.e. such
that $u+t\eta$ is convex for sufficiently small $t\ge 0$, we have,
$$\frac {d}{dt} A(u+t\eta)\big{|}_{t=0}
   =\frac {1}{n+2}\int_\Om w U^{ij} \p_{ij}\eta ,\tag 2.11$$
where $w=(\det \p^2 u )^{-(n+1)/(n+2)}$, and $[U^{ij}]$ is the
cofactor matrix of $\p^2 u$. Therefore a convex function is affine
maximal if and only if for any such $\eta$ with compact support in
$\Om$,
$$\int_\Om w U^{ij} \p_{ij}\eta \le 0 . \tag 2.12$$

We will need the upper semi-continuity of the affine area
functional. The upper semi-continuity was first proved in [15]. A
different proof was given in [21]. Here we give a simple proof.

\proclaim{Lemma 2.2} Let $\{u_m\}$ be a sequence of convex
functions in $\Om$, converging locally uniformly to $u$. Then
$$\limsup_{m\to\infty} A(u_m)\le A(u). \tag 2.13$$
\endproclaim

{\it Proof.} By the H\"older inequality we have
$$A(u, \Om)\le
 \bigg(\int \frac {\det \p^2 u}{\rho^{n+1}}\bigg)^{1/(n+2)}
 \bigg(\int \rho\bigg)^{(n+1)/(n+2)} \tag 2.14$$
for any positive function $\rho$. It follows
$$\align
A(u, \Om)=
 & \inf\bigg\{\bigg(\int \frac {\det \p^2 u}{\rho^{n+1}}\bigg)^{1/(n+2)},
 \ \ \rho\in C^0, \rho>0, \int\rho=1\bigg\}\\
 =& \inf\bigg\{\bigg(\int \frac {d\mu_r[u]}{\rho^{n+1}}\bigg)^{1/(n+2)},
 \ \ \rho\in C^0, \rho>0, \int\rho=1\bigg\} . \\
 \endalign $$
Since the singular part is defined on a set of measure zero, we
have
$$A(u, \Om)=
 \inf\bigg\{\bigg(\int \frac {d\mu[u]}{\rho^{n+1}}\bigg)^{1/(n+2)},
 \ \ \rho\in C^0, \rho>0, \int\rho=1\bigg\}  .\tag 2.15$$
The upper semi-continuity then follows from the weak continuity of
the Monge-Amp\`ere measure. $\square$

\vskip5pt

In the following we prove the existence and uniqueness of
maximizers for (1.8). Let $\Om$ be a bounded, Lipschitz domain in
$\R^n$, and $\phi$ a convex function (not necessarily smooth)
defined in a neighborhood of $\bom$. Denote by $\ol S[\phi, \Om]$
the set of convex functions $v$ satisfying $v=\phi$ on $\pom$ and
$N_v(\Om) \subset N_\phi(\bom)$. The latter relation means that if
we extend $v$ to the neighborhood of $\bom$ such that $v=\phi$
outside $\Om$, then $v$ is convex in the neighborhood of $\Om$.
Hence for any $\phi_1\in \ol S[\phi, \Om]$, we have $\ol S[\phi_1,
\Om]=\ol S[\phi, \Om]$ if we extend $\phi_1$ to a neighborhood of
$\Om$ by letting $\phi_1=\phi$.  If $\phi$ is uniformly convex,
then $\ol S[\phi, \Om]$ is the closure of $S[\phi]$ under uniform
convergence, where $S[\phi]$ is introduced in Section 1.

We will consider a more general maximization problem, that is
$$\sup\{A_f(v,\Om)\ {|} \ v\in \ol S[\phi, \Om]\} , \tag 2.17$$
where
$$A_f(v,\Om)=A(v,\Om)-\int_\Om fv \tag 2.18$$
and $f\in L^\infty(\Om)$ is a bounded, measurable function.

The existence of maximizers for (2.17) follows immediately from
the upper semi-continuity of the affine area functional. To prove
the uniqueness we first show that a maximizer has nonsingular
Monge-Amp\`ere measure.

\proclaim{Lemma 2.3} Let $u$ be a maximizer of (2.17). Then the
Monge-Amp\`ere measure $\mu[u]$ has no singular part.
\endproclaim

{\it Proof.} Suppose $\mu[u]$ has non-vanishing singular part
$\mu_s[u]$.  Since $\mu_s[u]$ is supported on a set of measure
zero and $\mu_r[u]$ is an integral function, it follows that for
any positive constant $K\ge 1$, there is a ball $B_r\subset\Om$
such that
$$\mu_s[u](B_r)\ge K\mu_r[u](B_r)+ 2 K^2|B_r| ,$$
for otherwise by the covering lemma, $\mu_s$ must be an integrable
function. Let $v$ be the solution to the Dirichlet problem
$$\align
\mu[v] & =K\mu_r[u]+ 2K^2\ \ \ \text{in}\ \ B_r,\\
v & =u\ \ \ \ \text{on}\ \ \p B_r .\\
\endalign $$
The existence and uniqueness of generalized solutions of the above
problem are well known, see [17]. Since $\mu[u] (B_r)\ge
K\mu_r[u](B_r)+ 2K^2 |B_r|$,  the set $E=\{v>u\}$ is not empty.
Let $\wtt u=u$ in $\Om-E$ and $\wtt u=v$ in $E$. Then $\wtt u\in
\ol S[\phi, \Om]$. We have
$$\align
A_f(v, E)-A_f(u, E)
 & = \int_E (\det \p^2 v)^{1/(n+2)}-\int_E (\det \p^2 u)^{1/(n+2)}
    -\int_E f(v-u)\\
 &\ge \int_E (K\det\p^2 u + 2K^2)^{1/(n+2)}
     -\int_E (\det \p^2 u)^{1/(n+2)} - C|E| .\\
\endalign $$
It is easy to see that the right hand side is positive, by
considering respectively the sets $\{\det \p^2 u>K\}$ and $\{\det
\p^2 u<K\}$, and choosing $K$ sufficiently large. Hence we obtain
$$A_f(\wtt u, \Om)=A_f(u, \Om-E)+A_f(v, E) > A_f(u, \Om) .$$
It follows that $u$ is not affine maximal, a contradiction.
$\square$

\proclaim{Theorem 2.1} Let $\Om$ be a bounded, Lipschitz domain in
$\R^n$. Suppose $\phi$ is a convex Lipschitz function defined in a
neighborhood of $\bom$, and $f$ is a bounded measurable function.
Then there is a unique maximizer $u$ for (2.17)
\endproclaim

{\it Proof.} As remarked above, the existence follows from the
upper semi-continuity of the affine area functional. To see the
uniqueness we observe that by the concavity of the affine area
functional, if both $u$ and $v$ are maximizers we have $\det \p^2
u=\det \p^2 v$ almost everywhere, and so $\mu[u]=\mu[v]$ as both
of them have no singular part by Lemma 2.3. It follows $u=v$ by
the uniqueness of generalized solutions to the Monge-Amp\`ere
equation [17]. We remark that the proof for the uniqueness of
generalized solutions in [17] (see Theorem 5.1 there) does not use
the strict convexity of the domain. $\square$

\vskip10pt

\noo{\it Remark 2.1}. Theorem 2.1 can be extended to the
functional
$$A_f(v, \Om)=A(v,\Om)-\int_\Om f(x, v),\tag 2.19$$
where $f(x, t)$ is locally bounded in $\bom\times \R^1$,
measurable in $x$ and convex in $t$. By the uniqueness in Theorem
2.1, one also sees that if $u_k$, $k=1, 2, \cdots$, are maximizers
of $\sup\{A_{f_k}(v,\Om)
           \ {|} \ v\in \ol S[\phi_k, \Om]\},$
and if $\phi_k\to \phi$ in $C^1(\bom)$, $f_k\to f$ in
$L^\infty(\Om)$, then $u_k\to u$ and $u$ is the maximizer of
(2.17).

\vskip30pt

\centerline{\bf\S 3. The general Plateau problem}

\vskip10pt

To study the maximization problem (1.3), we need to deal with
non-smooth, locally convex hypersurfaces. By definition, a locally
convex hypersurface is the image of a locally convex immersion in
$\R^{n+1}$ of a connected manifold $\N$, that is $\M=T(\N)$, on
which there is a continuous vector field pointing everywhere to
the convex side. This latter condition rules out hypersurfaces
such as $x_{n+1} =x_1\max(|x_1|-1, 0)$. Recall that an immersion
$T: \N\to\R^{n+1}$ is called locally convex if for any point
$p\in\N$, there is a neighborhood $N_\delta(p)\subset \N$ such
that $T(N_\delta(p))$ is a convex graph in $\R^{n+1}$. We say a
locally convex hypersurface $\M$ is convex if $\M$ lies in the
boundary of its convex closure.

For any given point $x$ on a locally convex hypersurface $\M$,
$T^{-1}(x)$ may contain more than one point in $\N$. To avoid
confusion in the following, when referring to a point $x\in\M$, we
need to understand a pair $(x, p)$, where $p=p_x\in \N$ such that
$T(p)=x$. Also we say $\omega_x\subset\M$ is a neighborhood of $x$
if it is the image of a neighborhood of $p$ in $\N$. We say
$\gamma$ is a curve on $\M$ if it is the image of a curve in $\N$
and so on. The $r$-neighborhood of $x$, $\omega_r(x)$, is the
connected component of $\M\cap B_r(x)$ containing the point $x$.

As a prelude we proved in [22] a fundamental result for locally
convex hypersurfaces, which plays a crucial role in our
investigation of the affine Plateau problem.

\proclaim {Lemma 3.1}    Let $\M$ be a compact, locally convex
hypersurface in $\R^{n+1}$, $n>1$. Suppose the boundary $\p\M$
lies in the hyperplane $\{x_{n+1}=0\}$. Then any connected
component of $\M\cap\{x_{n+1}<0\}$ is convex.
\endproclaim

Lemma 3.1 has two important applications to the affine Plateau
problem, namely to the existence of maximizers for (1.3) and their
strict convexity. Indeed, by Lemma 3.1 we conclude that the set
$S[\M_0]$ is precompact, and so by the upper semi-continuity of
the affine area functional, (1.3) admits a maximizer in $\ol
S[\M_0]$,  the closure of $S[\M_0]$ under local uniform
convergence. For the strict convexity of the maximizer, Lemma 3.1
enables us to reduce consideration to the graph case (see the
proof of Theorem 8.3).

To see that the set $S[\M_0]$ is precompact, we need a uniform
cone property of locally convex hypersurfaces.  Let $\C_{x, \xi,
r, \alpha}$ denote the cone with vertex $x$, axis $\xi$, radius
$r$, and aperture $\alpha$, that is,
$$\C_{x, \xi, r, \alpha}
 =\{y\in\R^{n+1}\ {|}\ |y-x|<r,\ \lan y-x, \xi\ran\ge \cos
 \alpha \, |y-x| \}. $$
We say that $\C_{x, \xi, r, \alpha}$ is an inner contact cone of
$\M$ at $x$ if this cone lies on the concave side of
$\omega_r(x)$.  We say $\M$ satisfies the {\it uniform cone
condition} with radius $r$ and aperture $\alpha$ if $\M$ has an
inner contact cone at all points with the same $r$ and $\alpha$.
From Lemma 3.1, we have

\proclaim{Lemma 3.2} Let $\M\subset B_R(0)$ be a locally convex
hypersurface with boundary $\p\M$. Suppose $\M$ can be extended to
$\wtt \M$ such that $\p\M$ is embedded in $\wtt\M$ and $\wtt\M-\M$
is locally strictly convex.  Then there exist $r, \alpha>0$
depending only on $n$, $R$, and the extended part $\wtt\M-\M$,
such that the $r$-neighborhood $\omega_r(x)$ is convex for any
$x\in\M$, and $\M$ satisfies the uniform cone condition with
radius $r$ and aperture $\alpha$.
\endproclaim

Lemma 3.2 was also proved in [22]. The main point of Lemma 3.2 is
that $r$ and $\alpha$ depend only on $n, R$ and the extended part
$\wtt\M-\M$. Therefore it holds with the same $r$ and $\alpha$ for
a family of locally convex hypersurfaces, which includes all
locally uniformly convex hypersurfaces, contained in $B_R(0)$,
with boundary $\p\M$ and Gauss mapping image coinciding with that
of $\M$. For any sequence of locally convex hypersurfaces in this
family, the uniform cone property implies that the sequence
sub-converges and the limit hypersurface is locally a graph. This
property was crucial for our resolution of the Plateau problem for
prescribed constant Gauss curvature [22] and also plays a key role
in the following existence proof of maximizers to the affine
Plateau problem.

\vskip10pt

Let $\M_0\subset \R^{n+1}$ be a bounded hypersurface with smooth
boundary which is smooth and locally uniformly convex up to its
boundary $\Ga$. As in [22] we extend $\M_0$ to a smooth, locally
uniformly convex hypersurface $\wtt\M_0$ such that $\Ga$ lies in
the interior of $\wtt\M_0$. Denote $\M_0^c=\wtt\M_0-\M_0$.

As in Section 1 we denote by $S[\M_0]$ the set of locally
uniformly convex hypersurfaces $\M$ with boundary $\Ga$, which can
be smoothly deformed from $\M_0$ in the family of locally
uniformly convex hypersurfaces whose Gauss mapping images lie in
that of $\M_0$. The latter assumptions is equivalent to saying
that $\M\cup \M_0^c$ is a locally convex hypersurface.  Let $\ol
S[\M_0]$ be the closure of $S[\M_0]$ under local uniform
convergence. By Lemma 3.2, $\ol S[\M_0]$ is well defined.
Moreover, any locally convex hypersurface in $\ol S[\M_0]$
satisfies the uniform cone condition with $r$ and $\alpha$
depending only on $n, R$ and $\M_0^c$.  From Lemma 3.1, we have
the following diameter estimate.

\proclaim{Lemma 3.3} If the image of the Gauss mapping of $\M_0$
does not cover any hemi-sphere in $S^n$, then there is $R>0$ such
that $\M\subset B_R(0)$ for any $\M\in S[\M_0]$ and so also for
any $\M\in \ol S[\M_0]$.
\endproclaim

For the extended affine Plateau problem, we are concerned with the
existence and regularity of maximizers to the problem
$$\sup_{\M\in \ol S[\M_0]} A(\M) .\tag 3.1$$

\proclaim{Theorem 3.1} Let $\M_0\subset \R^{n+1}$ be a bounded
hypersurface with smooth boundary which is smooth and locally
uniformly convex up to its boundary $\Ga$. Suppose the image of
the Gauss mapping of $\M_0$ does not cover any hemi-sphere in
$S^n$. Then there is a locally convex maximizer to (3.1).
\endproclaim

\noo{\it Proof}. \ The Gauss mapping image of $\M_0$, $\N$, is a
locally uniformly convex hypersurface immersed in $S^n$, and can
be decomposed into $m$ pieces, $\N=\bigcup^m_{i=1} F_i$, such that
every piece is strictly contained in some hemi-sphere, namely
every $F_i$ is a graph with uniformly bounded gradient. For any
$\M\in S[\M_0]$, let $\M^{(i)}=G^{-1}(F_i)$, where $G$ is the
Gauss mapping of $\M$. Then $\M=\bigcup \M^{(i)}$ and $\M^{(i)}$
is a graph for any $i$.

It follows by Lemma 3.3 that any sequence in $S[\M_0]$ contains a
convergent subsequence. Indeed, let $\{\M_j\}$ be a sequence in
$S[\M_0]$. For any given $j$, as above we decompose $\M_j$ into
the union of $\M_j^{(i)}$, $i=1, \cdots, m$. For each fixed $i$,
we can choose a coordinate system such that $\M_j^{(i)}$ are
graphs for all $j$. Hence by Lemma 3.3 and convexity, $\M_j^{(i)}$
contains a convergent subsequence. Moreover, by Lemma 3.3,
$A(\M_j^{(i)})$ is uniformly bounded for each $i$. Hence $A(\M_j)$
is uniformly bounded.

By the mollification of convex functions (see (2.3)) and the above
decomposition, it is easy to see that
$$\sup_{\M\in\ol S[\M_0]}A(\M)
     =\sup_{\M\in S[\M_0]}A(\M)<\infty.\tag 3.2$$
Hence the existence of maximizers to (3.1) follows from the upper
semi-continuity of the affine area functional. $\square$

The necessity of the condition in Theorem A, that is if a
hypersurface $\M$ is affine maximal, then the image of its Gauss
mapping cannot contain any hemi-sphere, is readily shown. Indeed,
if the Gauss mapping image contains the south hemi-sphere, we
denote by $\M'$ the preimage of the south hemi-sphere, given as a
graph of a convex function $u$ over a domain $\Om$. Then
necessarily $\det D^2 u\to\infty$ on $\pom$ and so $w\to 0$ on
$\pom$. Applying the maximum principle to equation (1.5),
regarding it as a linear, second order elliptic equation in $w$,
we find that $w\equiv 0$ in $\Om$. This is impossible. One can
also easily show that if the Gauss map image of $\M_0$ contains a
hemi-sphere, then the supremum in (3.1) is unbounded.

The rest of the paper is devoted to the regularity of maximizers
in Theorems 2.1 and 3.1.

\vskip30pt

\centerline{\bf \S 4. A priori estimates for classical solutions}

\vskip10pt

If the maximizer in Theorem 2.1 is smooth and locally uniformly
convex, it satisfies the nonlinear fourth order partial
differential equation
$$L[u]=f,\tag 4.1$$
where $L$ is the operator given in (1.5).

In this section we establish a priori estimates for solutions of
(4.1). These estimates are essentially proved in [21] where the
case $f\equiv 0$ is considered.

\proclaim{Lemma 4.1}   Let $u\in C^4(\Om)\cap C^{0,1}(\bom)$ be a
locally uniformly convex solution of (4.1) with $u=0$ on $\pom$.
Then, for any point $y\in\Om$, we have the estimate
$$\det D^2 u(y)\le C, \tag 4.2$$
where $C$ depends on $n$, $\dist (y, \pom)$, $\sup_\Om |Du|$,
$\sup_\Om f$, and $\sup_\Om |u|$.
\endproclaim

\noo{\it Proof.}\ \ Lemma 4.1 is proved in [21] (for the case
$f\equiv 0$). We include the proof here for completeness.

Let
$$z=\log \frac w{(-u)^\beta} -A|Du|^2, $$
where $\beta$ and $A$ are positive constants to be determined.
Then $z$ attains its minimum at an interior point $x_0\in\Om$. At
$x_0$ we have
$$\align
& 0=z_i=\frac {w_i}w -\beta \frac {u_i}u -2Au_ku_{ki},\\
& 0\le z_{ij}=\frac {w_{ij}}w-\frac {w_iw_j}{w^2}-\beta\frac
{u_{ij}}u +\beta\frac {u_iu_j}{u^2}-2Au_{ki}u_{kj}-2Au_ku_{kij}\\
\endalign $$
as a matrix. From $z_i=0$ we have
$$\frac {w_iw_j}{w^2} = \beta^2\frac {u_iu_j}{u^2}
  +\frac {2\beta A}{u} (u_iu_ku_{kj}+u_ju_ku_{ki})
  +4A^2u_ku_lu_{ki}u_{lj}. $$
Using the identities
$$\align
& u^{ij} u_{kij} = -\frac {1}{1-\th}\frac {w_k}w,
                                \ \ \ \th=\frac {1}{n+2} \\
& u^{ij}w_{ij}=U^{ij}w_{ij}/d =f/d,\ \ \ d=\det D^2 u,\\
\endalign $$
where $u^{ij}=U^{ij}/d$ is the inverse matrix of $D^2 u$, we
obtain
$$\align
 0 & \le u^{ij}z_{ij} \\
& = \frac {f}{d^\th} -\frac {\beta n}{u} -
    \frac {u^{ij}w_iw_j}{w^2} + \frac {\beta u^{ij} u_iu_j}{u^2}
    -2Au^{ij}u_{ki}u_{kj}+\frac {2A}{1-\th}\frac {u_k w_k}{w}\\
& = \frac {f}{d^\th} -\frac {\beta n}{u}
    -\beta(\beta-1)\frac {u^{ij} u_iu_j}{u^2}
    -2A\Delta u +\frac {4A^2 \th}{1-\th} u_{ij}u_iu_j
    -2\beta A \frac {1-2\th}{1-\th} \frac {|Du|^2}{u} \\
& \le \frac {f}{d^\th} -A\Delta u -\frac {\beta n}{u}
    +2\beta A \frac {|Du|^2}{|u|}, \\
\endalign$$
with the choice
$$A=\frac {1-\th}{4\th \sup_\Om |Du|^2} . $$
Therefore we have
$$ |u| \Delta u \le C(1+|Du|^2).$$
Hence $z(x_0)\ge -C$ if we choose $\beta\ge n(n+1)/(n+2)$. It
follows that $z(x)\ge z(x_0)\ge -C$ and so (4.2) holds. $\square$

\noo{\it Remark 4.1}. If $n=2$, the assumption $u=0$ on $\pom$ in
Lemma 4.1 can be removed. Indeed, let
$$z=\log \frac {w}{\eta^\beta}-A|Du|^2, $$
where $\eta(x)=(r^2-|x|^2)$ is a cut-off function and $r>0$ is
such that $B_r(0)\subset \Om$. Then similarly as above we have, at
a maximum point of $z$,
$$0\le u^{ij}z_{ij}
   \le \frac {f}{d^\th} -A\Delta u +\frac {\beta u^{ii}}{\eta}
    +2\beta A \frac {u_i\eta_i}{\eta}. $$
Note that
$$u^{11}+u^{22}=\frac {u_{11}+u_{22}}{u_{11}u_{22}}
              =\frac {\Delta u}{d}. $$
Hence we also obtain (4.2).

\noo{\it Remark 4.2.} Lemma 4.1 holds if $f=f(x, u, Du, D^2u)$ and
$f$ satisfies
$$f(x, z, p, r)\le C(1+\text{tr}\, r)\tag 4.3$$
for any symmetric matrix $r$.

\proclaim{Lemma 4.2}
 Let $u\in C^4(\Om)\cap C^{0,1}(\bom)$ be a locally uniformly convex
solution of (4.1). Suppose there exists an open set $\omega
\subset\Om$ such that $x\cdot Du<u$ in $\omega$ and $x\cdot Du= u$
on $\p\omega$. Then for any $y\in\omega$,
$$\det D^2 u(y)\ge C , \tag 4.4$$
where $C>0$ depends on $n$, $\dist (y, \p\omega)$, $\sup_\Om
|Du|$, $\inf_\Om f$ and $\sup_\omega |u-x\cdot D u|$.
\endproclaim

{\it Proof}.\ \ Let
$$z=\log w +\beta \log (u-x\cdot Du)+A|x|^2$$
for some positive constants $\beta$ and $A$ to be determined.
Suppose $z$ attains its maximum at $x_0\in\omega$. Then at $x_0$,
$$\align
& 0=z_i=\frac {w_i}w -\beta \frac {x_ku_{ki}}\phi + 2Ax_i,\\
& 0\ge z_{ii}=\frac {w_{ii}}w-\frac {w_i^2}{w^2}-\beta\frac
 {x_ku_{kii}+u_{ii}}\phi -\beta\frac {x_i^2u_{ii}^2}{\phi^2}+2A,\\
\endalign $$
where $\phi=u-x_iu_i$. By a rotation of coordinates we may suppose
$D^2 u$ is diagonal at $x_0$. Then
$$\align
0 \ge & u^{ii}z_{ii}\\
 = & \frac {f}{d^\th} - u^{ii} [\beta^2 \frac {x_i^2 u_{ii}^2}{\phi^2}
  -4\beta A \frac {x_i^2 u_{ii}}\phi + 4A^2 x_i^2]
  - \frac {\beta x_k}{\phi} u^{ii}u_{kii} -\frac {\beta n}{\phi}
  - \frac { \beta x_i^2 u_{ii}}{\phi^2} +2Au^{ii}\\
 = & \frac {f}{d^\th} -\beta(\beta+1)\frac {x_i^2 u_{ii}}{\phi^2}
  + \frac {4\beta A x_i^2}{\phi} + 2A (1-2A x_i^2) u^{ii}
  -\frac {\beta n}{\phi} + \frac {\beta x_i} {\phi (1-\th)}
     \frac {w_i}{w}\\
 = & \frac {f}{d^\th} - \frac {\beta n}\phi
    + \frac {4\beta A x_i^2}{\phi} + Au^{ii}
    - \beta (\beta+1) \frac {x_i^2 u_{ii}}{\phi^2}
    +\frac {\beta x_i}{\phi(1-\th)} (\beta \frac {x_i u_{ii}}\phi
          -2A x_i)\\
 \ge & \frac {f}{d^\th} - \frac C\phi +Au^{ii}
    +\beta(\frac \beta {1-\th} -\beta-1) \frac {x_i^2 u_{ii}}{\phi^2}\\
 \ge & \frac {f}{d^\th} - \frac C\phi +Au^{ii}\\
\endalign $$
if $\beta$ is large and $A$ is sufficiently small. It follows that
$|\phi| u^{ii} \le C$. Hence Lemma 4.2 holds. $\square$

We note that Lemma 4.2 also follows from Lemma 4.1, using the
Legendre transform, as in [21]; (see equation (7.4)). We can
determine neighborhoods $\omega=\omega_y$ of points $y\in\Om$
verifying the hypothesis of Lemma 4.2 in terms of a modulus of
convexity of the strictly convex function $u$. For any $y\in \Om$,
$h>0$, we define the section $S^0_{h, u}(y)$ by
$$S^0_{h, u}(y)=\{x\in\Om\ {|}\  u(x)<u(y)+Du(y)(x-y)+h\}.$$
The {\it modulus of convexity} of $u$ is a nonnegative function,
defined by
$$\rho_u (r) = \inf_{y\in\Om} \rho_{u, y}(r),
                                       \ \ \ \ r>0,\tag 4.5$$
where
$$\rho_{u, y}(r) = \sup
 \{h\ge 0\ \big{|} \ S^0_{h,u}(y)\subset B_r(y)\}  $$
if  there exists $h\ge 0$ such that $S^0_{h,u}(y)\subset B_r(y)$,
otherwise we define $\rho_{u, y}(r) = 0$. A general convex
function $u$ is strictly convex if and only if $\rho_u(r)>0$ for
all $r>0$.

If $u$ is a strictly convex solution of (4.1), we can characterize
the open set $\omega_y$ ($y\in\Om$) in the following way. Let
$\eps>0$ be any given constant. Let $\Cal P_\eps$ denote the set
of linear functions $g$ such that $g<u$ in $\Om$ and
$g(y)=u(y)-\eps$. Let $\ol g(x)=\sup\{g(x)\ |\ g\in\Cal P_\eps\}$.
Then $\ol g\le u$ and the graph of $\ol g$ is a convex cone. Let
$\omega$ denote the component of $\{\ol g<u\}$ containing $y$.
Then if $\eps<\rho_u(\frac 12 r)$, where $r=\dist(y, \pom)$, we
have $\ol \omega\subset \Om$.

By Lemmas 4.1 and 4.2 we have the following H\"older and Sobolev
space estimates.

\proclaim{Theorem 4.1} ({\it $W^{4, p}$ estimate)}\ \ Let $u\in
C^4(\Om)$ be a locally uniformly convex solution of (4.1). Then
for any $\Om'\subset\subset \Om$, $p\ge 1$, we have the estimate
$$\|u\|_{W^{4, p}(\Om')}\le C ,\tag 4.6$$
where $C$ depends on $n, p, \sup_\Om |f|$, $\dist(\Om', \pom)$,
and the modulus of convexity of $u$.
\endproclaim

\proclaim{Theorem 4.2} ({\it Schauder estimate})\ \ Let $u\in
C^4(\Om)$ be a locally uniformly convex solution of (4.1) with
$f\in C^\alpha(\bom)$, $0<\alpha<1$. Then $u\in C^{4,
\alpha}(\Om)$ and for any $\Om'\subset\subset \Om$,
$$\|u\|_{C^{4, \alpha}(\Om')}\le C ,\tag 4.7$$
where $C$ depends on $n, \alpha, \|f\|_{C^\alpha(\Om)}$,
$\dist(\Om', \pom)$, and the modulus of convexity of $u$.
\endproclaim

To prove Theorems 4.1 and 4.2, we have, by Lemmas 4.1 and 4.2 and
our control of the strict convexity of $u$,
$$C_1\le \det D^2 u\le C_2\tag 4.8$$
in any $\Om'\subset\subset\Om$. Now we write (4.1) as a second
order partial differential system
$$\align
U^{ij}w_{ij} & =f \ \ \ \text{in}\ \ \Om,\tag 4.9\\
\det D^2 u & =w^{-(n+2)/(n+1)}\ \ \ \text{in}\ \ \Om . \tag 4.10\\
\endalign$$
In [4] the authors proved, for the case $f\equiv 0$, a H\"older
estimate, with respect to sections of $u$, for solutions $w$ of
(4.9), assuming that the Monge-Amp\`ere measure $\mu[u]$ satisfies
a continuity condition with respect to the Lebesgue measure, which
is guaranteed by (4.8). By examining their argument one sees that
under (4.8), their H\"older continuity result holds when $f\in
L^\infty$. Taking account of the modulus of convexity of $u$, we
thus conclude the H\"older estimate for the function $w$. By the
interior Schauder estimate of the Monge-Amp\`ere equation [3], we
then obtain interior H\"older estimates for the second derivatives
of $u$. Hence (4.9) becomes a linear uniformly elliptic equation
with H\"older continuous coefficients. It follows we can estimate
$w$ in $W^{2, p}_{loc}(\Om)$ for any $p<\infty$. The $W^{4,
p}_{loc}(\Om)$ and $C^{4, \alpha}(\Om)$ estimates for $u$ now
follow from standard elliptic regularity theory.

\noo{\it Remark 4.3}.\ Theorems 4.1 and 4.2 hold for any strictly
convex solution of (4.1), by our regularity result, Theorem 6.2.

To prove the regularity of maximizers in Theorems 2.1 and 3.1, it
suffices to prove, in view of Theorems 4.1 and 4.2, that (a) the
maximizers are strictly convex and (b) they can be approximated by
smooth maximizers. We will prove (a) for dimension $n=2$ in
Section 8 and (b) for all dimensions in Section 6.

\vskip30pt

\centerline{\bf \S 5. The second boundary value problem}

\vskip10pt

In this section we prove the existence of solutions to the
following boundary value problem,
$$\align
L[u] & =f (x, u)  \ \ \ \  \text{in}\ \ \Om ,  \tag 5.1 \\
   u & =\phi\ \ \ \  \text{on}\ \ \pom , \tag 5.2 \\
   w & =\psi\ \ \ \  \text{on}\ \ \pom , \tag 5.3 \\
\endalign$$
where $L$ and $w$ are as in (1.5) (1.6), $\Om$ is a smooth,
uniformly convex domain in $\R^n$, $\phi, \psi$ are smooth
functions on $\pom$ with
$$C_0^{-1} \le \psi \le C_0      \tag 5.4$$
for some positive constant $C_0$. We suppose $f\in
L^\infty(\Om\times\R)$, $f$ is non-decreasing in $u$, and there is
$t_0\le 0$ such that
$$ f(x, t)\le 0 \ \ \ \text{when}\ \ t\le t_0.\tag 5.5$$

\proclaim{Theorem 5.1}\ \ The boundary value problem (5.1)-(5.3)
admits a solution $u\in W^{4,p}_{loc}(\Om)$ $\cap C^{0, 1}(\bom)$
($\forall\ p>1$) with $\det D^2 u\in C^0(\bom)$. If $f\in
C^\alpha(\bom\times \R)$, where $\alpha\in (0, 1)$, then the
solution $u\in C^{4,\alpha}(\Om)\cap C^{0, 1}(\bom)$.
\endproclaim

We will use Theorem 5.1 in the next section to construct smooth
approximations to the maximizers in Theorem 2.1.

To prove Theorem 5.1 we write (5.1) as a system (4.9) (4.10) and
consider the approximating problem
$$\align
U^{ij} w_{ij} & =f\ \ \ \ \ \ \text{in}\ \ \Om, \tag 5.6\\
\det D^2 u & =\eta w^{-(n+2)/(n+1)}+(1-\eta)
                           \ \ \text{on}\ \ \Om,\tag 5.7\\
\endalign $$
where $u$ and $w$ satisfy the boundary condition (5.2) and (5.3),
and $\eta=\eta_k\in C_0^\infty(\Om)$ is a nonnegative cut-off
function satisfying $\eta=1$ in $\Om_k=\{x\in\Om\ {|}\ \dist(x,
\pom)<1/k\}$.

\proclaim{Lemma 5.1} Let $(u, w)$ be a $C^2$ smooth solution of
(5.6)-(5.7). Then there exists a constant $C>0$ such that
$$\align
 C^{-1} \le w & \le C\ \ \ \text{in}\ \ \Om, \tag 5.8\\
 |w(x)-w(x_0)|& \le C|x-x_0|
        \ \ \ \forall\ \ x\in\Om, x_0\in\pom, \tag 5.9\\
\endalign$$
where $C$ depends only on $n$, $\diam(\Om)$, $\sup_\Om |f|$, and
$\sup_\Om |u|$, and is independent of $k$.
\endproclaim

\noo{\it Proof.}  Let $z=\log w -h(u)$, where $h$ is a convex,
monotone increasing function satisfying $h(t)=t$ when $t\ge -t_0$
and $h(t)\ge -t_0-1$ when $t\le -t_0$. If $z$ attains its minimum
at a boundary point, by (5.4) we have $w\ge C$ in $\Om$. If $z$
attains its minimum at an interior point $x_0\in\Om$. At this
point we have
$$\align
0 & = z_i=\frac {w_i}{w}-h'(u) u_i,\\
0 &\le z_{ij} = \frac {w_{ij}}{w}-\frac {w_iw_j}{w^2}
              -h'(u)u_{ij}-h''(u)u_iu_j\\
\endalign$$
as a matrix. Hence
$$0\le u^{ij} z_{ij} \le \frac {f}{d^\th} -nh'(u)$$
where $d=\det D^2 u$, $\th=1/(n+2)$. If $u(x_0)\le t_0$, $f=0$ and
we reach a contradiction. Hence $u(x_0)\ge t_0$ and $h'(u)>0$. We
obtain $d(x_0)\le C$. Since $z(x)\ge z(x_0)$, we obtain
$$w(x)\ge w(x_0) \text{exp} (h(u(x))-h(u(x_0))). $$
The first inequality in (5.8) follows.

Next let $z=\log w+A|x|^2$. If $z$ attains its maximum at a
boundary point, by (5.4) we have $w\le C$ and so (5.8) holds. If
$z$ attains its maximum at an interior point $x_0$, we have, at
$x_0$,
$$\align
0 & = z_i=\frac {w_i}{w}+2Ax_i,\\
0 &\ge z_{ii} = \frac {w_{ii}}{w}-\frac {w_i^2}{w^2}+2A .\\
\endalign$$
Suppose $(D^2 u)$ is diagonal at $x_0$. Then
$$\align
0 & \ge u^{ij} z_{ij}\\
 & = \frac {f}{dw} -4A^2 x_i^2 u^{ii}+2Au^{ii}\\
 &\ge \frac {f}{dw } +Au^{ii} \\
 \endalign $$
if $A$ is small. Hence
$$dw\sum u^{ii}\le C.$$
By (5.7) we obtain
$$[\eta w^{-(n+2)/(n+1)}+(1-\eta)]^{(n-1)/n}w
           \le dw\sum u^{ii}\le C.\tag 5.10$$
We obtain $w\le C$, and hence (5.8) is proved.

Let $v$ be a smooth, uniformly convex function in $\Om$ such that
$v=\psi$ on $\pom$ and $D^2 v\ge K$. Then
$$U^{ij}v_{ij} \ge K \sum U^{ii}\ge CK[\det D^2 v]^{(n-1)/n}\ge CK.$$
Hence if $K$ is large enough, $v$ is a lower barrier of $w$ by
applying the comparison principle to (5.6).  We thus obtain
$$w(x)-w(x_0)\ge -C|x-x_0|
        \ \ \ \forall\ \ x\in\Om, x_0\in\pom.\tag 5.11$$
Similarly one can construct an upper barrier for $w$. Hence (5.9)
holds. $\square$

By approximation, Lemma 5.1 holds for $w\in W^{2, p}(\Om)$ with
$p>n$. Indeed, let $\{f_k\}$ be a sequence of bounded functions
which converges to $f$ in $L^p$, and let $w_k$ be the solution of
(5.6) with $f=f_k$, where $U^{ij}$ is the cofactor matrix of $D^2
u$, which is independent of $k$. Then $w_k\to w$ in $W^{2, p}$. As
above we have the estimate $dw_k\sum u^{ii}\le C$. Sending
$k\to\infty$ we obtain (5.10) and so the second inequality in
(5.8) follows. The first inequality in (5.8) can be proved in the
same way as above.

\proclaim{Lemma 5.2} There is a solution $(u, w)$, where  $u\in
C^{2, \alpha}(\bom)$ and $w\in W^{2, p}(\Om)$ ($p>n$), to the
above approximation problem.
\endproclaim

\noo{\it Proof.}\ By (5.8), $u$ is uniformly bounded, and is
strictly convex in $\Om$ [2]. Applying the interior H\"older
continuity result [4] to (5.6) we have $\det D^2 u\in
C^\alpha(\Om)$ for some $\alpha\in (0, 1)$, which in turn implies
$w\in W^{2, p}(\Om)$ $\forall\ p>1$ and $u\in C^{2, \alpha}(\Om)$.
Near the boundary we also have $u\in C^{2, \alpha}$ by applying
the regularity theory of Monge-Amp\`ere equation to (5.7)
[5,10,13]. Therefore we have global regularity for the
approximation problem. Next we use the degree theory to prove the
existence of solutions.

For any positive $w\in C^{0, 1}(\bom)$, let $u=u_w$ be the
solution of (5.7) with $u=\phi$ on $\pom$, and let $w_t$, $t\in
[0, 1]$, be the solution of
$$\align
 U^{ij}w_{ij} & = tf(x, u)\ \ \text{in}\ \ \Om ,\tag 5.12\\
          w_t & =t\psi+(1-t)\ \ \ \text{on}\ \ \pom .\\
\endalign $$
Then the mapping $T_t:\ w\in C^{0,1}(\bom) \to w_t\in
C^{0,1}(\bom)$ is compact. By the above a priori estimates, the
degree $\text{deg} (T_t, B_R, 0)$ is well defined, where $B_R$ is
the set of all positive function satisfying
$\|w\|_{C^{0,1}(\bom)}\le R$. When $t=0$, from (5.12) we have
obviously $w\equiv 1$. Namely $T_0$ has a unique fixed point
$w\equiv 1$. Hence the degree $\text{deg} (T_t, B_R, 0)=1$ for all
$t\in [0,1]$. This completes the proof. $\square$

Denote by $u_k$ (corresponding to $\eta_k$) the solution obtained
in Lemma 5.2, where the estimates for the upper and lower bounds
of $\det D^2 u_k$ are independent of $k$. Letting $u_k\to u$ we
obtain (5.8) for $w=[\det D^2 u]^{-(n+1)/(n+2)}$.  Therefore we
conclude $u\in W^{4, p}_{loc}(\Om)\cap C^{0,1}(\bom)$ by Theorem
4.1. Note that by (5.8) we also have (5.9), which implies that
$w\in C^0(\bom)$ and $u$ satisfies the boundary conditions (5.2)
and (5.3). If $f\in C^\alpha(\bom\times\R)$, then $u\in C^{4,
\alpha}(\Om)\cap C^{0, 1}(\bom)$ by Theorem 4.2. Hence we obtain
Theorem 5.1.

In a separate paper [23] we prove the uniqueness and global
regularity of solutions to the boundary value problem (5.1)-(5.3).

\vskip30pt

\centerline{\bf \S 6. Approximation by smooth solutions}

\vskip10pt

In this section we show that an affine maximal function can be
approximated by smooth solutions of the affine maximal surface
equation. Our approach also embraces the inhomogeneous case and
utilizes the solvability of the second  boundary value problem
(Theorem 5.1).

We begin by considering a particular version of the second
boundary value problem. Let $\phi\in C^2(\ol B) $ be a uniformly
convex function in a ball $B=B_R(0)$, vanishing on $\p B$. Let
$H\in C^\infty ( - \infty, 1)$ be a non-negative convex function
such that
$$ H(t) = \cases
4^{n}(1-t)^{-2n} \ \ \ &\text{if}\ \ 1/2<t<1,\\
t^4 &\text{if}\ \  t<-1.\\
\endcases \tag 6.1$$
Let $f\in L^\infty(\Om)$ be the function in (2.18) and suppose
$\Om\subset\subset B$. Extend $f$ to $B$ such that
$$f(x, u)= h(u-\phi(x))\ \ \ \text{in}\ \ B-\Om, \tag 6.2$$
where $h(t)=H'(t)$. Then $f$ is nondecreasing in $u$.

\proclaim {Lemma 6.1}   There is a locally uniformly convex
solution $u$ to the second boundary value problem
$$\align
L[u] & =f(x, u)\ \ \ \text{in}\ \ B, \tag 6.3\\
  u  & = \phi \ \ \ \text{on}\ \ \p B,\\
  w  & =1 \ \ \ \text{on}\ \ \p B\\
\endalign $$
with $u\in W^{4, p}_{loc}(B)\cap C^{0,1}(\ol B)$, for all
$p<\infty$, $w\in C^0(\bom)$, where $L$ is the operator in (1.5).
\endproclaim

\noo{\it Proof}. We will prove that if $u$ is a locally uniformly
convex solution of (6.3), then
$$\inf_B u\ge -K_0\tag 6.4$$
for some $K_0>0$ depending only on $\phi$ and $R$, the radius of
the ball $B$; and
$$|f(x, u)|\le C . \tag 6.5$$
Once (6.4) and (6.5) are established, Lemma 6.1 follows from
Theorem 5.1.

First we prove (6.4). Let $\delta>0$ be a small constant. Since
$\phi$ vanishes on $\p B$, $\Om_\delta=\{u<-\delta\}$ is strictly
contained in $B$. We compute
$$\align
\int_{\Om_\delta} U^{ij}w_{ij}(u+\delta)
 & = -\int_{\Om_\delta} (u+\delta)_i w_j U^{ij} \tag 6.6\\
 & = -\int_{\pom_\delta} (u+\delta)_i\gamma_j w U^{ij}
     +\int_{\Om_\delta} w U^{ij} (u+\delta)_{ij}\\
 & < n\int_B w\det D^2 u =n\int_B [\det D^2 u]^\th, \\
\endalign $$
where $\th=1/(n+2)$, $\gamma$ is the unit outward normal, and we
have used the divergence free property of $[U^{ij}]$ for any fixed
$i$ or $j$. Sending $\delta$ to 0,  we obtain
$$\int_B f(x, u)u\le n\int_B [\det D^2 u]^\th
 = n\int_\M K^\th , $$
where $\M$ is the graph of $u$, and $K$ is the Gauss curvature of
$\M$. It follows that
$$\int_B f(x, u)u
  \le n|\M|^{1-\th} \big[\int_\M K\big]^\th
  \le C|\inf u|^{1-\th} . $$
Recalling that $f$ is bounded in $\Om$, we obtain
$$\int_{B-\Om} f(x, u) u \le C|\inf u|. $$
Since $u$ is convex and $u=0$ on $\p B$, there exists $C>0$ such
that
$$|\inf_{B-\Om} u|^4\ge C\int_B |u|^4. $$
It follows by our construction of $H$, see (6.1),
$$\int_{B-\Om} f(x, u) u\ge C_1 |\inf u|^4-C_2$$
for some positive constants $C_1, C_2>0$. Hence (6.4) holds.

Next we prove (6.5). Since $f(x, t)$ is increasing in $t$, $f(x,
u)$ is bounded from below by (6.4). If suffices to prove that
$f(x, u)$ is bounded from above.

We first prove $\det D^2 u$ is bounded near $\p B$. Indeed,  by
convexity and our choice of $H$, $f(x, u)$ is bounded from above
near $\p B$. For any boundary point $x_0\in \p B$, we suppose by a
rotation of axes that $x_0=(R, 0, \cdots, 0)$. Let
$\ell(x)=ax_1+b$ be a linear function such that
$\ell(x_0)<u(x_0)=0$ and $\ell>u$ on $x_1=R-\delta_0$, where
$\delta_0>0$ is a constant such that $f$ is upper bounded in
$B\cap\{x_1>R-\delta_0\}$. Let $z=\log \frac {w}{u-\ell}$. If $z$
attains a minimum at a boundary point $\p B$, by the boundary
condition $w=1$ in (6.3) we see that $z$ is bounded from below and
so $\det D^2 u$ is bounded from above near $\pom$. If $z$ attains
a minimum at some interior point $y_0\in \{u>\ell\}$, we compute,
at $y_0$,
$$\align
0= & z_i=\frac {w_i}{w}-\frac {(u-\ell)_i}{u-\ell},\\
   & z_{ij}=\frac {w_{ij}}{w}-\frac {(u-\ell)_{ij}}{u-\ell}\\
\endalign $$
with the matrix $[z_{ij}]\ge 0$. It follows that,
$$0\le u^{ij}z_{ij}=\frac {f}{[\det D^2 u]^\th}-\frac {n}{u-\ell}. $$
Hence $\det D^2 u\le C$ at $y_0$ and so $z\ge C$, which in turn
implies that $\det D^2 u$ is bounded near $\p B$.

By (6.4) we then conclude that $u$ is Lipschitz at $\p B$, and
hence $Du$ is uniformly bounded in $B$ by convexity. Returning to
(6.6), we have
$$\align
\int_{\Om_\delta} U^{ij}w_{ij} (u+\delta)
      & \ge -\int_{\pom_\delta} u_i\gamma_j w U^{ij} \\
      & =  -\int_{\pom_\delta} u_\gamma w U^{\gamma\gamma}\\
      & = -\int_{\pom_\delta} (u_\gamma)^n w K_s\\
  & \ge - n\omega_n \big(\sup_B |Du|\big)^n\inf_{\pom_\delta} w\\
\endalign $$
where $K_s$ denotes the Gauss curvature of $\pom_\delta$. Letting
$\delta\to 0$, we obtain from (6.3),
$$\int_B (-f(x, u)) u\le C\tag 6.7$$
since $w=1$ on $\p B$.

If $u(x)-\phi(x)$ is sufficiently close to one at some point $x\in
B-\Om$, then $u(x)-\phi(x)$ is sufficiently close to one at nearby
points in $B-\Om$, by the convexity of $u$ and $\phi$. Hence the
integral on the left hand side of (6.7) must become very large,
which is in contradiction with the estimate (6.7). Hence (6.5)
holds. Lemma 6.1 now follows from Theorem 5.1. $\square$

We remark that the function $f(x, u)$ is not defined when $u\ge
\phi+1$. This is not a problem for the use of the degree argument
in the proof of Lemma 5.2. One can also choose a sequence $f_j(x,
u)$ which is defined for all $u\in\R$ and converges to $f$.

Next we use Lemma 6.1 and the penalty method to prove the
maximizer in Theorem 2.1 can be approximated locally by smooth
local maximizers.

\proclaim{Theorem 6.1} Let $\Om$ and $\phi$ be as in Theorem 2.1.
Then for any convex $\Om'\subset\subset\Om$, there exists a
sequence of smooth solutions of equation (4.1), ($\in W^{4,
p}_{loc}(\Om)\ \forall\ p<\infty$), converging uniformly in $\Om'$
to the maximizer $u$.
\endproclaim

\vskip5pt

\noo{\it Proof}. \ \ Without loss of generality let us assume that
$\Om$ is convex,  $u$ is a maximizer in $\Om^\delta = \{x\in\R^n\
{|}\ \dist(x, \Om)<\delta\}$ for some $\delta>0$ small, $\phi=u$
in $\Om^\delta-\Om$. We will prove that $u$ can be approximated by
smooth solutions of (4.1) in $\Om$.

Let $B_R=B_R(0)$ be a ball in $\R^n$ containing $\ol{\Om^\delta}$.
Let
$$\wtt\phi (x)=\sup_{v\in\Cal P} v(x),
                              \ \ \ \ x\in B_R , $$
where $\Cal P$ is the set of linear functions $v$ such that $v\le
\phi$ in $\Om^\delta$ and $v\le K_0$ in $B_R$ for some given
constant $K_0$. Since $\phi$ is Lipschitz and $\Om$ is convex, we
can choose $K_0$ sufficiently large such that $\wtt \phi=\phi$ in
$\Om^\delta$. By definition, $\wtt\phi$ cannot be strictly convex
at any point in $B_R-\ol{\Om^\delta}$. By Aleksandrov's theorem,
the set $N_{\wtt\phi}(B_R-\ol{\Om^\delta})$ has measure zero, that
is $\mu[\wtt\phi]=0$ in $B_R-\ol{\Om^\delta}$. Therefore we may
suppose directly that $\phi$ is defined in $B_R$ such that
$\mu[\phi]=0$ in $B_R-\ol{\Om^\delta}$ and $\phi$ equals the
constant $K_0$ on $\p B_R$.

Let $\{\phi_k\}$ be a sequence of convex functions such that
$\phi_k = \phi$ in $\ol{\Om^\delta}$, $\phi_k = \phi$ on $\p B_R$,
$\phi_k$ is uniformly convex in  $\ol B_R-\ol{\Om^\delta}$, and
$\phi_k \to \phi$ uniformly in $\ol B_R$. Let $H_j (t) = H(2^jt)$
be a sequence of smooth, convex functions, defined in $(-\infty,
2^{-j})$.  Let $f_{k, j}(x, u)= f(x)$ when $x\in\Om$ and $f_{k,
j}(x, u) = h_j(u-\phi_k)$ when $x\in B_R-\Om$, where $h_j=H'_j$
and $f$ is the function in (2.17).

By Lemma 6.1, there is a convex solution $u_{k, j}$ of (6.3) with
$f=f_{k, j}$, which is an extremal of the concave functional
$$J_{k,j} (u)=J_{k, j}(u, B_R) = A(u, B_R)-\int_\Om f u
             - \int_{B_R-\Om}  H_j(u-\phi_k),$$
where $A(u, B_R)$ is the affine area functional on the domain
$B_R$.

Similar to (6.4) we have $u_{k,j}\ge -K_0$ for some $K_0$
independent of $k, j$. We have indeed a stronger estimate, for any
given $k$,
$$\inf_{B_R-\Om} (u_{k, j}-\phi_k)\to 0
     \ \ \ \text{as}\ j\to\infty.\tag 6.8$$
To prove (6.8) we suppose $\inf_{B_R-\Om} (u_{k, j}-\phi_k)$ is
attained at $x_j$. Let $\ell$ be the tangent plane of $\phi_k$ at
$x_j$. Let $\omega=\{u_{k,j}<\ell\}$. We compute
$$\align
\int_\omega U^{ij}w_{ij}(u-\ell)
 & = -\int_\omega (u-\ell)_iw_j U^{ij}\ \ \ \ (u=u_{k,j})\\
 & = -\int_{\p\omega} (u-\ell)_i\gamma_j w U^{ij}
     + \int_\omega (u-\ell)_{ij}w U^{ij}\\
\endalign $$
The first integral on the right hand side is negative. Hence we
obtain
$$\align
\int_\omega f_{k,j}(x, u) (u-\ell)
 & =\int_\omega U^{ij}w_{ij}(u-\ell)\\
 & \le \int_\omega w\det D^2 u \\
 & = \int_\omega (\det D^2 u)^\th\\
 & =\int_\M K^\th \le C,\\
 \endalign $$
where $\M$ is the graph of $u$ and $K$ is the Gauss curvature. If
(6.8) is not true, the integral on the left hand side converges to
infinity (as $j\to\infty$) by our definition of $f_{k, j}$, which
is a contradiction. Hence (6.8) holds.

Observe that by our definition of $H_j$,
$$u_{k, j}\le \phi_k+2^{-j} ,\tag 6.9 $$
which, together with (6.8), implies that
$$u_{k, j}\to\phi_k \tag 6.10$$
in $B_R-\Om$ as $j\to\infty$.

Since $\phi_k\to\phi$, we have $u_{k, j}\to \phi$ ($k,
j\to\infty$) locally uniformly in $B_R-\Om$ as long as $j$ is
large enough. By convexity, $u_{k,j}$ sub-converges to a convex
function $\ol u$ in $B_R$. By our definition of $\ol S[\phi,\Om]$
in Section 2, the function $\ol u$, when restricted to $\Om$,
belongs to the set $\ol S[\phi, \Om]$. We want to prove that $\ol
u$ is the maximizer of (2.17), whence $\ol u=u$ by the uniqueness
assertion in Theorem 2.1.

Choose $R'<R$ and $r>0$ such that $\Om\subset B_{R'-3r}(0)$.
Denote
$$D_\sigma=\{x\in B_R(0)\ |\ R'-r-\sigma<|x|<R'+r+\sigma\}$$
where $0\le \sigma\le r$. Let
$$\ol u_{k, j}=\sup\{v\ |\ v\in\Phi_{k, j}\},$$
where $\Phi_{k, j}$ is the set of convex functions in $B_R(0)$
which satisfy $v\le u_{k, j}$ in $B_R-D_\sigma$, $v\le\max(\phi_k,
u_{k, j})$ in $D_\sigma$, and $v\le \phi_k$ in $D_0=D_{\sigma\,|
\sigma=0}$. Then for any fixed $r>\sigma>0$, by the uniform
convexity of $\phi_k$ and (6.10), we have
$$\align
\ol u_{k, j} &=u_{k, j}\ \ \ \text{in}\ \ B_R-D_\sigma,\tag 6.11\\
\ol u_{k, j} &=\phi_k\ \ \ \text{in}\ \
    \{x\in D_0\ |\ \dist(x, \p D_0)>\sigma/2\}, \\
|\ol u_{k, j} & -\phi_k|\le |u_{k, j}-\phi_k|\ \ \ \text{in}\ \
D_\sigma. \tag 6.12\\
\endalign $$
provided $j$ is sufficiently large. By (6.11), $\ol u_{k, j}\in
\ol S[u_{k, j}, B_R]$. Since $u_{k, j}$ is the maximizer of $J_{k,
j}$  in $\ol S[u_{k, j}, B_R]$, we have
$$J_{k, j}(\ol u_{k, j})\le J_{k, j}(u_{k, j}), \tag 6.13$$
Observe that
$$\align
\int_{D_\sigma} [\det D^2 u]^{1/(n+2)}
 & \le |D_\sigma|^{(n+1)/(n+2)}
      \big{|}\int_{B_R} \det D^2 u\big{|}^{1/(n+2)}\\
 & \le C|D_\sigma|^{(n+1)/(n+2)},\\
 \endalign $$
where $u=u_{k, j}$. By (6.12) we have
$$J_{k, j}(\ol u_{k, j})\ge J_{k, j}(u_{k, j})-\eps\tag 6.14$$
with $\eps\to 0$ as $r, \sigma\to 0$.

Let
$$v_{k, j}=\cases
u\ \ \  &\text{in}\ \Om^\delta,\\
\phi_k\ \ \ &\text{in}\ \ B_{R'}-\Om^\delta,\\
\ol u_{k, j}\ \ &\text{in}\ \ B_R-B_{R'}\\
\endcases $$
be an extension of $u$ in $B_R$. Then $v_{k,j}=u=\phi$ in
$\Om^\delta-\Om$, and $v_{k, j}\in \ol S[u_{k, j}, B_R]$. Hence we
have, by (6.13), (6.14),
$$J_{k, j}(v_{k, j})\le J_{k, j}(u_{k, j})
       \le J_{k, j}(\ol u_{k, j})+\eps.  $$
Since $v_{k, j}=\ol u_{k, j}$ in $B_R-B_{R'}$, we obtain
$$J_{k, j}(v_{k, j}, B_{R'})
           \le J_{k, j}(\ol u_{k, j}, B_{R'})+\eps.  \tag 6.15$$
Since $H_{k, j}\ge 0$, we have
$$\align
J_{k,j}(\ol u_{k,j}, B_{R'})
   & \le A(\ol u_{k,j}, B_{R'})-\int_\Om f\ol u_{k,j}\\
   & = A_f(\ol u_{k,j}, \Om)+ A(\ol u_{k,j}, B_{R'}-\Om), \tag 6.16\\
\endalign $$
where
$$A_f(u, \Om)=A(u, \Om)-\int_\Om fu. $$
We compute
$$\align
A(\ol u_{k,j}, B_{R'}-\Om)
 & = A(\ol u_{k,j}, \Om^\delta-\Om)
       + A(\ol u_{k,j}, B_{R'}-\Om^\delta)\\
 & \le  |\Om^\delta-\Om|^{(n+1)/(n+2)}
  \big[\int_{\Om^\delta-\Om} \det D^2 \ol u_{k,j} \big]^{1/(n+2)}\\
   &\ \ \ + C\big[\int_{B_{R'}-\Om^\delta}
                 \det D^2\ol u_{k,j} \big]^{1/(n+2)} .\\
   \endalign $$
We choose $\delta>0$ sufficiently small such that the first term
on the right hand side is small. Observe that $\ol u_{k, j}\to
\phi$ and the Monge-Amp\`ere measure $\mu[\phi]=0$ in
$B_R-\Om^\delta$. We therefore obtain
$$ A(\ol u_{k,j}, B_{R'}-\Om) \le 2\eps \tag 6.17$$
if $k$ is large enough. Therefore by (6.15), (6.16), and the upper
semi-continuity of the affine area functional, we obtain
$$\align
A_f(u, \Om) &\le  J_{k, j}(v_{k, j}, B_{R'})\\
      &\le J_{k, j}(\ol u_{k, j}, B_{R'})+\eps\\
      &\le A_f(\ol u_{k,j}, \Om) +2\eps\\
      &\le A_f(\ol u, \Om)+3\eps \\
\endalign $$
if $k, j$ are sufficiently large. Hence $\ol u$ is the maximizer,
and thus $\ol u=u$. $\square$

The penalty method is used above to create a sequence of smooth
solutions of (4.1) satisfying (6.10). The function $H$, chosen in
(6.1), facilitates the estimate (6.9). This function can be
replaced by convex functions defined in $(-\infty, \infty)$ which
grow much faster for $t>0$ than for $t<0$.

From Theorem 6.1 we obtain the following interior regularity in
all dimensions, which includes the case of affine maximal
functions and thereby completes the proof of Theorem B.

\proclaim{Theorem 6.2}  Let $u$ be a strictly convex maximizer of
(2.17), with $f\in C^\infty(\Om)$. Then $u\in C^\infty (\Om)$ and
satisfies equation (4.1) in $\Om$.
\endproclaim

In the case of an affine maximal graph $\M=\M_u$, (that is
$f\equiv 0$), we obtain from Theorem 6.1 and the uniqueness part
of Theorem 2.1, that there exists a sequence of smooth functions
$\{u_m\}\subset C^\infty(\Om)$ with affine maximal graphs,
converging locally uniformly to $u$. As a byproduct, we may extend
our affine Bernstein result in [21] to non-smooth affine maximal
surfaces.

\proclaim{Theorem 6.3} Let $u$ be an affine maximal convex
function defined in the whole space $\R^2$. Suppose $u$ has at
least one strictly convex point. Then $u$ is a quadratic function.
\endproclaim

The assumption that $u$ has at least one strictly convex point
implies that the level set $S^0_{h, u}(p)$ is a bounded convex
domain for some point $p\in\R^2$. Note that if $u$ has no strictly
convex point, then $u(x_1, x_2)=\phi(x_1)$ for some convex
function $\phi$ in an appropriate coordinate system.

In the following two sections we will show that the hypothesis of
strict convexity can be dispensed with in two dimensions.

\vskip30pt

\centerline{\bf \S 7. The generalized Legendre transform}

\vskip10pt

In this section we introduce the Legendre transform for
(nonsmooth, locally) convex functions in general domains, which
will be used in the next section to prove the strict convexity of
maximizers in Theorems 2.1 and 3.1.

Let $\Om$ be a convex domain and $u\in C^2(\bom)$ be a uniformly
convex function. The Legendre transform of $u$ is the function
$u^*$ defined in $\Om^*=Du(\Om)$, given by
$$\align
u^*(y) & =\sup\{\wtt x\cdot y-u(\wtt x)\ |\ \wtt x\in\Om\}\\
       & = x\cdot y-u(x), \tag 7.1\\
       \endalign $$
where $x\in\Om$ is uniquely determined by $y=Du(x)$. The Legendre
transform $u^*$ is a uniformly convex, $C^2$ smooth function in
$\Om^*$. Furthermore the Legendre transform of $u^*$ is $u$
itself.

From the relation $y=Du(x)$ we have $x=D u^*(y)$ and
$$\{D^2 u(x)\}=\{D^2 u^*(y)\}^{-1}. \tag 7.2$$
Therefore if $u\in C^2(\bom)$ is affine maximal, namely if it is a
maximizer of the affine area functional $A$ in the set $\ol S[u,
\Om]$, $u^*$ is a maximizer of the associated functional
$$A^*(v)=\int_{\Om^*} [\det \p^2 v]^{(n+1)/(n+2)} \tag 7.3$$
in the set $\ol S[u^*, \Om^*]$. If $u$ is a smooth, uniformly
convex solution of the non-homogeneous equation (4.1), then by
direct computation, $u^*$ satisfies the equation
$$U^{ij} w^*_{ij}=-\frac {1}{n+1} f(Du^*)\,\det D^2 u^*, \tag 7.4$$
where $U^{ij}$ is the cofactor matrix of $D^2 u^*$ and $w^*=[\det
D^2 u^*]^{-1/(n+2)}$.

Next we extend the Legendre transform to non-smooth (locally)
convex functions.

Let $\Om$ be a bounded $C^2$ smooth domain in $\R^n$ and $\phi$ be
a (locally) uniformly convex function in $\bom$. As before we
denote by $S[\phi, \Om]$ the set of locally uniformly convex
functions $u\in C^2(\Om)\cap C^{0, 1}(\bom)$ satisfying $u=\phi$
and $Du(\Om)\subset D\phi(\Om)$, and denote by $\ol S[\phi, \Om]$
the closure of $S[\phi, \Om]$ under local uniform convergence. In
the following we extend the Legendre transform to all functions in
$\ol S[\phi, \Om]$ and show that if $u$ is a maximizer for
$\sup_{v\in \ol S[\phi, \Om]} A(v)$, the Legendre transform of $u$
is a maximizer for $\sup_{v\in \ol S[\phi^*, \Om^*]} A^*(v)$,
where $\Om^*=D\phi(\Om)$.

First we consider the case when the domain $\Om$ is convex such
that $D\phi$ is a diffeomorphism from $\bom$ to $D\phi(\bom)$.
Extend $\phi$ to a neighborhood of $\bom$, say, the domain
$\Om^\sigma=\{x\in\R^n\ {|}\ \dist(x, \Om)<\sigma\}$ for some
$\sigma>0$, such that $\phi$ is uniformly convex in $\Om^\sigma$
(see [22]). For any $u\in \ol S[\phi, \Om]$, we extend $u$ to
$\Om^\sigma$ such that $u=\phi$ in $\Om^\sigma-\Om$. Let
$$u^*(y)=\sup_{x\in\Om^\sigma} (x\cdot y -u (x))
             \ \ \ \forall\ y\in D\phi(\Om^\sigma) , \tag 7.5$$
where the supremum is attained at a point $x$ such that $y\in
N_u(x)$. It is easy to see that $u^*=\phi^*$ in
$D\phi(\Om^\sigma-\Om)$. We define the Legendre transform of $u$
as the restriction of $u^*$ on the domain $\Om^*=D\phi(\Om)$.

Let $\ol S[\phi^*] = \ol S[\phi^*, \Om^*]$ denote the set of
convex functions $v$ such that $v=\phi^*$ on $\pom^*$ and
$N_v(\Om^*)\subset \bom$. For any convex function $v\in \ol
S[\phi^*, \Om^*]$, we extend $v$ to $D\phi(\Om^\sigma)$ such that
$v=\phi^*$ in $D\phi(\Om^\sigma-\Om)$. Then similarly we can
define the Legendre transform of $v$, which is a convex function
on $\Om$.

\noo{\it Remark 7.1.} When $\Om$ is convex and $\phi$ is a convex
function defined in $\Om^\sigma$, the Legendre transform (7.5) is
well defined for any $u\in\ol S[\phi, \Om]$. The smoothness of
$\pom$ and $\phi$ is not required. This implies that if $u$ is a
convex function in a non-convex domain $\Om$ and $Du(x)=Du(y)$ for
two points $x,y\in\Om$ such that the line segment $\ol{xy}$ is
contained in $\Om$, then the tangent planes of $u$ at $x$ and $y$
coincide.

We remark that for any function $u\in S[\phi, \Om]$, the usual
Legendre transform of $u$ is defined by
$$\hat u^*(y)=\sup_{x\in\Om^\sigma} (x\cdot y -u (x))
             \ \ \ \forall\ y\in Du(\Om) .$$
The function $\hat u^*$ is a convex function defined in
$Du(\Om)\subset\Om^*$. Obviously we have $u^*=\hat u^*$ in
$Du(\Om)$. The graph of $u^*$ in $\Om^*-Du(\Om)$ consists of line
segments. In fact we have
$$u^*(x)=\sup \ell(x),$$
where the supremum is taken over all linear function satisfying
$\ell\le \hat u$ in $Du(\Om)$ and $\ell\le \phi^*$ on $\pom^*$.
Our definition of the Legendre transform is an extension of the
usual one such that for any function $u\in\ol S[\phi, \Om]$, $u^*$
is properly defined in the whole $\Om^*$.

\vskip10pt

Next we consider the case when $\Om$ is a general bounded smooth
domain. In this case the gradient mapping $D\phi$ may not be one
to one. We need to regard $D\phi$ as an immersion and the domain
$\Om^*=D\phi(\Om)$ as an immersed manifold in $\R^n$. Let
$S^0[\phi, \Om]$ denote the set of uniformly convex functions
$u\in S[\phi, \Om]\cap C^2(\bom)$ satisfying $u=\phi$ on $\pom$
and $Du(\Om)=D\phi(\Om)$. For any $u\in S^0[\phi, \Om]$, extend
$u$ to $\Om^\sigma$ such that $u=\phi$ outside $\Om$. Then $u$ has
continuous gradient across $\pom$.  It is easy to see that $\ol
S[\phi, \Om]$ is also the closure of $S^0[\phi, \Om]$ under
uniform convergence.

We claim that for any $u_0, u_1\in S^0[\phi, \Om]$, there exists a
diffeomorphism $\rho$ from $\Om$ into itself, such that
$Du_0(x)=Du_1(\rho(x))$. Indeed, for any $t\in [0, 1]$, the
function $u_t:=tu_1+(1-t)u_0\in S^0[\phi,\Om]$. Since $u_0$ and
$u_1$ are uniformly convex, namely they are $C^2$ up to boundary
and $D^2 u_0$ and $D^2 u_1$ are positive definite, there exists a
diffeomorphism $\rho_t$ from $\Om$ into itself, such that
$Du_t(x)=Du_{t'}(\rho_t(x))$ if $t<t'$ and $t, t'$ are
sufficiently close. Hence the claim follows.

It follows that $\Om^*$, as a manifold immersed in $\R^n$, is also
the image of the immersion $D u:\ \Om\to\R^n$ for any $u\in
S^0[\phi, \Om]$. For $u\in S^0[\phi, \Om]$, the Legendre transform
of $u$, $u^*$, is a {\it single valued} function defined on the
manifold $\Om^*$, such that for any $y\in \Om^*$,
$$u^*(y)=x\cdot y -u (x),\tag 7.6$$
where $x\in\Om$ is the unique point such that $y=Du(x)$. Obviously
(7.6) is an extension of (7.1).

\noo{\it Remark 7.2.}\ \ We need to explain how $u^*$ is
understood as a single valued function on $\Om^*$. Let us
introduce a manifold $\N=\{(x, y)\ |\ x\in\Om, y=Du(x)\}$ with
metric $ds^2=dy^2$. Then the gradient mapping $Du$ is an embedding
of $\Om$ to $\N$ and $Du(\Om)$ is isometric to $\N$. The function
$u^*$ can be regarded as a single valued function on $\N$.

Next we define the Legendre transform for (non-smooth) functions
$u\in \ol S[\phi, \Om]$. Let $\{u_k\}\subset S^0[\phi, \Om]$ be a
sequence of smooth, uniformly convex functions converging to $u$.
We define $u^*$, the Legendre transform of $u$, by
$$u^*(y)=\lim_{k\to\infty} u^*_k(y) ,\tag 7.7 $$
where $u^*_k$ is the Legendre transform of $u_k$, given in (7.6).
To show that $u^*$ is well defined, we need to show that $u^*$ is
independent of the choice of the sequences. Suppose
$\{u^1_k\}\subset S^0[\phi, \Om]$ is another sequence of smooth,
uniformly convex functions which converges to $u$. For any point
$y$ on the manifold $\Om^*$, there is a unique point $x_k\in\Om$
($x^1_k\in\Om$, resp.) such that $Du_k(x_k)=y$ ($Du^1_k(x^1_k)=y$,
resp.). By choosing subsequences we suppose $x_k\to x_0$ and
$x^1_k\to x^1_0$.

We claim that there exists a curve $\ell\subset\bom$ connecting
$x_0$ and $x^1_0$ such that $y\in N_u(x)$ for any $x\in\ell$.
Suppose for a moment the claim is true. Since $u$ can be extended
to a convex function in a neighborhood of $\Om$ (namely $u=\phi$
on $\Om^\sigma-\Om$), we can divide $\ell$ into the union of small
arcs, $\ell=\cup_i\ell_i$, such that for each $\ell_i$, the line
segment $\hat \ell_i$ connecting the two endpoints of $\ell_i$
lies in $\Om^\sigma$. Let $x_0, \cdots, x_m$ be the endpoints of
these line segments such that $\hat\ell_i=\overline{x_ix_{i+1}}$.
Since $y\in N_u(x_i)$, the hyperplane $P_i=\{(x,
x_{n+1})\in\R^{n+1}\ |\ x_{n+1}=y\cdot (x-x_i)+u(x_i)\}$ is a
tangent plane of $u$ at $x_i$. By convexity $u$ is linear on these
line segments $\hat \ell_i$ and the tangent planes of $u$ at $x_i$
and $x_{i+1}$ coincide, see Remark 7.1. Therefore we obtain
$x_0\cdot y -u (x_0)=x_0^1\cdot y -u (x_0^1)$, namely the limit in
(7.7) does not depend on the choice of the sequence $\{u_k\}$. By
Remark 7.1, one also sees that the Legendre transform depends
locally on $u$.

To prove the claim we let $u_k^t=tu_k+(1-t)u^1_k$. Then for any
$t\in [0, 1]$, $\{u^t_k\}\subset S^0[\phi, \Om]$ and $u^t_k\to u$
as $k\to\infty$. Let $x_k^t$ be the unique point in $\Om$ such
that $Du^t_k(x^t_k)=y$. By the uniform convexity of $u_k$ and
$u^1_k$, $x^t_k$ depends continuously on $t$. Hence for fixed $k$,
$E_k=\{x^t_k\ {|}\ \ t\in [0,1]\}$ is a curve in $\Om$. Let $\Cal
E$ denote the set of such points $\hat x$ that there is a sequence
$\{\hat x_k\}$, where $\hat x_k\in E_k$, which sub-converges to
$\hat x$. For any given $\hat x\in\Cal E$, since the hyperplane
$P_k=\{x_{n+1}=y\cdot (x- x^t_k)+u^t_k (x^t_k)\}$ is a tangent
plane of $u^t_k$ at $x^t_k$ and $u^t_k$ converges to $u$
uniformly, $\Cal P=\{x_{n+1}=y\cdot (x- \hat x)+u (\hat x)\}$ is a
tangent plane of $u$ at $\hat x$. It follows that $y\in N_u(\hat
x)$ for any $\hat x\in \Cal E$. Since $E_k$ is connected, for any
$r>0$ small such that $B_r(\hat x)\subset\Om^\sigma$ and $E_k\cap
\p B_r(\hat x)\ne \emptyset$, there is a point $\hat x_r\in \p
B_r(\hat x)\cap\Cal E$. Namely $\Cal P_r=\{x_{n+1}=y\cdot (x- \hat
x_r)+u (\hat x_r)\}$ is a tangent plane of $u$ at $\hat x_r$. When
$r$ is sufficiently small, by convexity we see that $\Cal P$ and
$\Cal P_r$ must coincide. Again by the convexity of $u$, we see
that $u$ is linear on the line segment $\ol{\hat x\hat x_r}$ and
so $y\in N_u(x)$ for any $x$ on the line segment. The claim is
proved. Hence (7.7) is well defined.

Regarding $\Om^*$ as a manifold (see Remark 7.2), we can define
the Legendre transform for functions $v\in \ol S[\phi^*, \Om^*]$
similarly. Indeed,  note that $S^0[\phi^*, \Om^*]$ is the set of
Legendre transforms of functions in $S^0[\phi, \Om]$, see (7.6).
Hence if $v\in S^0[\phi^*, \Om^*]$, we can define the Legendre
transform for a function $v\in S^0[\phi^*, \Om^*]$ by
$$v^*(x)=x\cdot y -v (y),\tag 7.8$$
where $y\in \Om^*$ is the unique point such that $x=Dv(y)$.
Obviously $v^*\in S^0[\phi, \Om]$ and $v$ is the Legendre
transform of $v^*$. For general functions $v\in \ol S[\phi^*,
\Om^*]$, the Legendre transform of $v$ is also defined as the
limit of the Legendre transforms of sequences of smooth, locally
uniformly convex functions in $S^0[\phi^*, \Om^*]$ which converge
to $v$. From (7.7), this is well defined.

\noo{\it Remark 7.3}.\ It is easy to see that if $\{u_k\}$ is a
sequence in $\ol S[\phi, \Om]$ and $u_k\to u$ uniformly, then the
Legendre transform $u^*_k$ converges to $u^*$. Observe that $u\in
S^0[\phi, \Om]$ if and only if $u^*\in S^0[\phi^*, \Om^*]$. Hence
for any $u\in \ol S[\phi, \Om]$, we have $(u^*)^*=u$, namely the
Legendre transform of $u^*$ is $u$ itself.

\vskip10pt

\proclaim{Lemma 7.1} Let $u\in\ol S[\phi, \Om]$ and let $u^*$ be
its Legendre transform. Then
$$ A(u, \Om) =A^*(u^*, \Om^*), \tag 7.9$$
where $\Om^*$ is regarded as a manifold, as noted in Remark 7.2.
\endproclaim

\noo{\it Proof}. If $u$ is twice differentiable at $x$ and $\det
\p^2 u(x)>0$, then it is easy to prove that $u^*$ is twice
differentiable at $y=Du(x)$, and by (7.2) one has the relation
$$\big(\p^2 u^*(y)\big) \big( \p^2 u(x)\big)=I , \tag 7.10$$
where $I$ is the unit matrix.

By our definition of the Legendre transform, (7.9) is obviously
true for uniformly convex functions in $S^0[\phi, \Om]$. We need
to prove (7.9) for general convex functions in $\ol S[\phi, \Om]$.
For any Borel set $E\subset \Om$ with $\dist (E, \pom)\ge \delta$,
and any $\alpha\in (0, 1)$, we have
$$\align
\int_E \big(\det \p^2 u\big)^\alpha
 &\le |E|^{1-\alpha}\bigg(\int_E\det \p^2 u\bigg)^\alpha \\
 & =  |E|^{1-\alpha} (\mu_r[u](E))^\alpha \\
 & \le |E|^{1-\alpha} (\mu[u](E))^\alpha \\
 &\le  |E|^{1-\alpha}
  \bigg(\frac {\text{osc}\, u}{\delta}\bigg)^{n\alpha},\\
\endalign$$
where $\mu_r[u]$ is the regular part of the Monge-Amp\`ere measure
$\mu[u]$. It follows that $(\det \p^2 u)^\alpha$ is locally
equi-integrable whenever $u$ is bounded.

Extend $u$ to $\Om^\sigma$ such that $u=\phi$ outside $\Om$. Let
$u_h$ be the mollification of $u$, as defined in (2.3). Then $D^2
u_h(x)\to \p^2 u(x)$ whenever $u$ is twice differentiable at $x$
[25]. Let $\hat u_h=u_h+h|x|^2$. Then we have
$$\int_\Om (\det \p^2 u)^\alpha
  =\lim_{h\to 0} \int_\Om (\det D^2  \hat u_h)^\alpha.\tag 7.11$$
Let $\hat u_h^*$ be the Legendre transformation of $\hat u_h$. By
(7.10) we have $D^2 \hat u_h^*(y)\to \p^2 u^*(y)$ whenever $u^*$
is twice differentiable at $y$ and $\det \p^2 u^*(y)>0$. It
follows
$$\int_{\Om^*} (\det \p^2 u^*)^\alpha
  \le \lim_{h\to 0} \int_{\Om^*} (\det D^2  \hat u_h^*)^\alpha.
                           \tag 7.12$$
Since (7.9) holds for uniformly convex functions in $S^0[\phi,
\Om]$, we obtain
$$A^*(u^*, \Om^*)\le A(u, \Om) . $$

For any $v\in \ol S[\phi^*, \Om^*]$, let $v^*$ be the Legendre
transform of $v$. Then similarly we have
$$A(v^*, \Om)\le A^*(v, \Om^*). $$
Observe that the Legendre transform of $u^*$ is $u$ itself.
Choosing $v=u^*$ such that $v^*=u$, we obtain (7.9). $\square$

\vskip5pt

Similar to Theorem 2.1, the functional $A^*$ is upper
semi-continuous, so there exists a maximizer $v$ to the supremum
$$\sup_{v\in \ol S[\phi^*,\Om^*]} A^*(v). \tag 7.13$$
If $v$ is smooth and uniformly convex, it satisfies the Euler
equation (7.4) with $f=0$.

Let $u$ be the Legendre transform of $v$. From (7.9) we see that
$u$ is the maximizer for $\sup_{u\in\ol S[\phi,\Om]} A(u)$, and
$$\sup_{u\in \ol S[\phi,\Om]} A(u)
      =\sup_{v\in \ol S[\phi^*,\Om^*]} A^*(v). \tag 7.14$$
Furthermore,  $v$ is the Legendre transform of $u$, and hence by
Theorem 2.1, $v$ is unique.

\vskip10pt

For locally convex hypersurfaces, the notion of support function
plays a similar role to the Legendre transform for convex graphs.
As a locally convex hypersurface can be decomposed as a union of
graphs, (see the proof of Theorem 3.1), we will just discuss
support functions for convex hypersurfaces, that is hypersurfaces
lying in the boundaries of convex bodies.

Let $\M$ be a smooth, convex hypersurface with Gauss mapping image
$\N$. The support function $H$ is a function defined on $\N$,
given by
$$H(x) = \sup \{ p \cdot x\ |\   p \in \M \}. \tag 7.15 $$
If $\M$ is locally uniformly convex, the supremum is attained at
the unique point $p\in\M$ with unit outward normal $x$. Moreover,
the Gauss curvature of $\M$ at $p$ is given by
$$K=1/\det (\D^2 H+HI)(x), \tag 7.16$$
where $\D $ is the covariant derivative on $S^n$ (under a local
orthonormal frame), and $I$ is the unit matrix. Extend $H$ to
$\R^{n+1}$ such that it is homogeneous of degree one, namely
$H(tx)=tH(x)$ for any $t\ge 0$ and $x\in\N$. Then $\M$ can be
recovered from $H$ by
$$\M=\{D H (x)\in \R^{n+1}:\ |\ x\in \N\},\tag 7.17$$
where $D$ is the ordinary derivative in $\R^{n+1}$. If locally
$\M$ is represented as a graph of a convex function $u$, then
$$u^*(y)= H (y,-1) \ \ \   y\in \Om^*, \tag 7.18$$
is exactly the Legendre transform of $u$. By (7.16) we have
$$\det D^2 u^*(y)=(1+|y|^2)^{-(n+2)/2}\det (\D^2 H+HI).\tag 7.19 $$
Note that in the graph case, (7.17) is equivalent to the Legendre
transform for $u^*$.

For an arbitrary convex hypersurface $\M$, one can define the
support function using the generalized Gauss mapping. For any
interior point $p\in\M$, the generalized Gauss mapping is a
multi-valued mapping given by
$$G(p)=\{x\in S^n\ |\
    x\cdot \wtt p\le x\cdot p\ \ \forall \wtt p\in\M\}. \tag 7.20$$
Let $\N=\bigcup G(p)$, where the union is taken over all interior
points of $\M$. Then we can define the support function of $\M$
also by (7.15), and if $\M$ is a graph, the support function is
related to the Legendre transform by (7.18).

\vskip30pt

\centerline{\bf \S 8. Strict convexity}

\vskip10pt

Throughout this section we assume the dimension $n=2$.  In two
dimensions, the local strict convexity of affine maximal functions
which vanish on the boundary of a convex domain $\Om$ follows from
our treatment of the affine Bernstein problem in [21] and the
approximation in Section 6.

\proclaim {Theorem 8.1} Suppose that $u\in C^0(\bom)$ is an affine
maximal convex function in a domain $\Om\subset\R^2$,  vanishing
on the boundary $\pom$. Then $u$ is strictly convex in $\Om$ and
the modulus of convexity of $u$ can be estimated in term of $\Om$
and $\inf_\Om u$.
\endproclaim

We remark that Theorem 8.1 holds for affine maximal graphs with
boundary lying on a plane. Using the affine invariance of the
affine maximal surface equation, the proof of Theorem 8.1 reduces
to the fact that a convex function with bounded Monge-Amp\`ere
measure is differentiable at any point on its graph, not lying on
a line segment joining two boundary points. This is embodied in
the following lemma.

\proclaim {Lemma 8.1} Suppose $u$ is a nonnegative convex function
in a domain $\Om\ni \{0\}$, $\subset\R^2$, satisfying  $u>0$ on
$\pom$, $u(0)=0$ and $u(x_1, 0)\ge |x_1|$. Then the Monge-Amp\`ere
measure $\mu[u]$ cannot be a bounded function.
\endproclaim

\noo{\it Proof}. We outline the proof here, as it is essentially
contained in [21] Section 5. First by Lemma 2.3, the Radon-Nikodym
derivative of $\mu[u]$ is an integrable function.  If the set
$\{u=0\}$ is the single point $\{0\}$, let $G_\eps=\{u<\eps\}$.
Then the Lebesgue measure $|G_\eps| \le \eps\delta_\eps$, where
$\delta_\eps=\sup\{|x_2|\ {|}\  (x_1, x_2)\in G_\eps\}$. On the
other hand, the image of the normal mapping $N_u(G_\eps)$
satisfies $|N_u(G_\eps)|>C\eps/\delta_\eps$ for some $C>0$
independent of $\eps$. By definition we have
$\mu[u](G_\eps)=|N_u(G_\eps)|$. Hence $\mu[u]$ is unbounded near
the origin.

If the set $\{u=0\}$ is a line segment, by the assumption that
$u>0$ on $\pom$, we may suppose $\{u=0\}$ is the line segment
$\{(0, -t)\ {|}\ \ 0\le t\le t_0\}$ such that the origin is an
endpoint. Let $G_\eps=\{u<\ell_\eps\}$, where
$\ell_\eps(x)=\eps+\delta_\eps x_2$ is a linear function, and
$\delta_\eps>0$ is chosen such that $u(0, \eps/\delta_\eps)
=\ell_\eps(0, \eps/\delta_\eps) =2\eps$ and $u(0,
-\eps/\delta_\eps)=\ell_\eps(0, -\eps/\delta_\eps) =0$. By the
convexity and the assumption that $u(x_1, 0)\ge |x_1|$, we see
that $G_\eps\subset \{-\delta_\eps\le x_2\le C\delta_\eps\}$ for
some $C>0$ depending on $\sup |Du|$. Similarly as above we have
$|G_\eps| \le C\eps\delta_\eps$, and
$|N_u(G_\eps)|>C\eps/\delta_\eps$ for a different $C>0$. Hence
$\mu[u]$ is unbounded near the origin. $\square$

\vskip5pt

We say a point $x_0\in \p F$ is an extreme point of a convex set
$F\subset \R^n$ if there is a hyperplane $P$ such that
intersection $P\cap \ol F$ is the single point $\{x_0\}$. The main
result of this section is the following strict convexity in
dimension 2 for affine maximal graphs with general boundary
values.

\proclaim{Theorem 8.2} Let $\Om$ be a bounded $C^2$ smooth domain
in $\R^2$ and $\phi$ be a uniformly convex function in $\bom$.
Then the maximizer $u$ in Theorem 2.1 is strictly convex in $\Om$,
in the case of vanishing $f$.
\endproclaim

\noo{\it Proof}.\ If $u$ is not strictly convex, then the graph of
$u$, $\M_u$, contains a line segment $\ell$. Let $P$ be the
tangent plane of $\M_u$ containing $\ell$. Let $F$ be the
component of $P\cap\M_u$ which contains $\ell$ (note that
$P\cap\M_u$ may contain more than one component if $\Om$ is not
convex). Then $F$ is a convex set. We have two cases.

Case a): $F$ has an extreme point $x_0$ which is an interior point
of $\M_u$. Then there is a plane $P_0=\{v(x) =a\cdot x +b\}$, such
that $v(x_0) > u(x_0)$ and $v<u$ on $\pom$. Let $\Om'=\{x\in\Om\
{|}\ v(x)>u(x)\}$. Then $u$ is strictly convex in $\Om'$ by
Theorem 8.1. This is a contradiction since $x_0$ is an extreme
point of $F$.

Case b): All extreme points of $F$ are boundary points of $\Om$.
In this case we suppose for simplicity that $P=\{x_3=0\}$ and
$u\ge 0$ near $F$.  Since all extreme points of $F$ are boundary
points of $\Om$, there exists a line segment contained in $F$,
which we suppose is
$$\ell=\{te_1\ \big{|}\ -1\le t\le 1\},\tag 8.1$$
where $e_1=(1, 0)$, such that both endpoints of $\ell$ are
boundary points of $\M_u$. By the uniform convexity of $\phi$ we
see that $\ell$ is transversal to $\pom$ at the endpoints $\pm
e_1$, namely $|\lan e_1, \gamma\ran|\ge C>0$, for otherwise we
have $\frac {\p^2}{\p x_1^2}\phi=0$ at the endpoints of $\ell$,
where $\gamma$ is the unit outward normal to $\pom$. Since
$\phi=u\ge 0$ on $\pom$ and $u = 0$ at the endpoints $\pm e_1$, by
the smoothness of $\phi$ and $\pom$, we have $|u(x)|\le \frac 12 C
x_2^2$ for $x\in \pom$, close to the endpoints $\pm e_1$. By
convexity it follows that
$$u(x)\le \frac C2 x_2^2\tag 8.2$$
for $x\in \Om$, close to the segment $\ell$.

Let $u^*$ be the Legendre transform of $u$, as introduced in
Section 7. Then $u^*$ is the maximizer of the functional $A^*$
over the set $\ol S[\phi^*, \Om^*]$, see (7.14), and the origin is
an interior point of $\Om^*$, where $\Om^*=D\phi(\Om)$. By (8.2)
and since $u\ge 0$, we have $u^*(0)=0$ and
$$u^*(x_1, x_2)\ge \frac 1{2C} x_2^2 \tag 8.3$$
near the origin. Since the line segment $\ell$ in (8.1) is
contained in $F$, we have
$$u^*(x_1, x_2)\ge |x_1|  \tag 8.4$$
locally near the origin. In the following we prove that the
function $u^*$, which satisfies (8.3) and (8.4), cannot be a
maximizer of $A^*$ over the set $\ol S[\phi^*, \Om^*]$. Before
continuing our proof, we make the following remarks.

\vskip5pt

\noo{\it Remark 8.1.}\  By Lemma 8.1 the Radon-Nikodym derivative
of the Monge-Amp\`ere measure $\mu[u^*]$ is an unbounded function.
On the other hand if one can prove that $u^*$ can be approximated
by smooth solutions of (7.4) (with $f=0$), then $\mu[u^*]$ is
bounded by Lemma 4.1 or Remark 4.1. However at the moment we don't
know if $u^*$ can be approximated by smooth solutions. In our
proof of the approximation result in Section 6, we used the power
$\frac 1{n+2}$ and the concavity of the affine area functional
$A$, see (5.10) and (6.13).

\noo{\it Remark 8.2}. Case b) can occur if $\phi$ is not uniformly
convex. Indeed, if the graph of $\phi$ contains a line segment of
which both endpoints are boundary points, then the line segment is
on the graph of any function in $\ol S[\phi]$.

\noo{\it Remark 8.3}. In Section 7 the Legendre transform of $u$
is defined on $D\phi(\Om)$, which is regarded as an immersed
manifold in $\R^n$. Here in the proof of Theorem 8.2, we need only
to consider the piece of $\M_{u, \delta}=\{(x, u(x))\in\M_u\ |\
u(x)<\delta\}$ for small $\delta>0$, where (as in \S 7) $u$ is
extended to a neighborhood of $\Om$ such that $u=\phi$ outside
$\Om$. The Legendre transform of $u$ introduced in \S 7, when
restricted to the piece $\M_{u, \delta}$, coincides with
$$u^*(y)=\sup_{x\in\omega} (x\cdot y -u (x))
             \ \ \ \forall\ y\in \omega^*=N_u(\omega)  $$
since $\omega$ is convex, where $\omega=\omega_\delta$ is the
projection of $\M_{u, \delta}$ on $\{x_{n+1}=0\}$. Moreover, $u$
can be recovered from $u^*$ by the same formula, that is
$$u(x)=\sup_{y\in\omega^*} (x\cdot y -u^*(y))
             \ \ \ \forall\ x\in \omega. $$
In the following we suppose directly that the Legendre transform
of $u$ is given by the above formulae.

\noo{\it Remark 8.4}. We shall use the following basic property of
the Legendre transform. If $x_0\in\bom$ and $y_0\in N_u(x_0)$,
then $x_0\in N_{u^*}(y_0)$ and locally $u^*(y)\ge x_0\cdot
y-u(x_0)$, namely $x_0\cdot y-u(x_0)$ is a support plane of $u^*$
at $y_0$, where $N_u$ is the normal mapping introduced in Section
2. Since the endpoints $\pm e_1$ of $\ell$ are boundary points of
$\pom$, we have $\lim_{x_1\to 0}\frac 1{|x_1|} u^*(x_1, 0)=1$.

\vskip10pt

Returning to the proof of Theorem 8.2, we denote $g(t)=u^*(0, t)$.
Then $g(0)=0$ and by (8.3)
$$g(t)\ge t^2 .\tag 8.5$$
Note that we may suppose $C\le 1/2$ in (8.3) by replacing $u$ by
$u/2C$ in (8.2). Denote by $\Om_h=\{x\in\R^2\ |\ u^*(x)<h\}$ the
sub-level set of $u^*$, where $h>0$. We also denote by $t^+_h>0$
and $t^-_h<0$ the unique constants satisfying $g(t^\pm_h)=h$.
First we show that $\Om_h$ has $``$good$"$ shape.

\proclaim{Lemma 8.2} There is a positive constant $C>0$,
independent of $h$, such that
$$\Om_h \subset \{x\in\R^2\ |\ -h\le x_1\le h,
                    Ct^-_h\le x_2\le Ct^+_h\}. \tag 8.6$$
\endproclaim

\noo{\it Proof.} By (8.4) we have $\Om_h \subset \{-h\le x_1\le
h\}$. We need only to prove $\Om_h \subset \{Ct^-_h\le x_2\le
Ct^+_h\}$. By (8.4), this is true with $C=2$ if $u^*(x_1,
0)=|x_1|$. In general we will prove $\Om_h \subset \{x_2\le
Ct^+_h\}$ by restricting to the region $\Om_h\cap \{x_2\ge 0\}$.
For $h>0$ small, let
$$u^*_h(x_1, x_2)=h^{-1}u^*(hx_1, t_h^+x_2).$$
Then the graph of $u^*_h$ (the part in $\{x_3<1\}$) sub-converges,
as $h\to 0$, to a convex surface $\M^*_0$ such that the point $(0,
1, 1)\in \p\M^*_0$. Since the endpoints of $\ell$ are boundary
points of $\Om$, we have $\lim_{x_1\to 0}\frac 1{|x_1|}u^*(x_1,
0)=1$. Hence the line segments $\{x_3=|x_1|, x_1\in (-1, 1),
x_2=0\}$ lie on $\M^*_0$. On the other hand, since $u^*(x_1,
x_2)\ge |x_1|$, we see that $\M^*_0$ lie above the planes $x_3=\pm
x_1$. By convexity it follows that $\M^*_0\subset \{y_2<2\}$.
Hence $\Om_h \subset \{x_2\le Ct^+_h\}$ with $C\to 2$ as $h\to 0$.
Similarly one can prove $\Om_h \subset \{x_2\ge Ct^-_h\}$ with
$C\to 2$ as $h\to 0$. $\square$

In the above proof, a sequence of convex surfaces $\{\Ga_k\}$ is
said to converge to a surface $\Ga_\infty$ if for any given $R >
0$ and any $\delta>0$, there exists $k>1$ such that for any $j\ge
k$, $\Ga_j\cap B_R(0)\subset N_\delta(\Ga_\infty\cap B_R)$ and
$\Ga_\infty\cap B_R(0)\subset N_\delta(\Ga_j\cap B_R)$, where
$N_\delta$ denotes the $\delta$-neighborhood.

By the proof of Lemma 2.3, the singular part of the Monge-Amp\`ere
measure $\mu [u]$ vanishes. Hence
$$\lim_{t\to 0} g(t)/|t|=0. $$

\proclaim{Lemma 8.3} There is a constant $\beta>0$ such that
$$g(t)\le |t|^{1+\beta}  \tag 8.7$$
for small $t$.
\endproclaim

\noo{\it Proof.} Denote $\ol g(t)=g(t)+g(-t)$, where $t>0$. It
suffices to prove that $\ol g$ satisfies (8.7). We need only to
prove that there is a positive constant $\th>0$, independent of
$h$ such that
$$\ol g(\frac 12 t)\le (\frac 12-\th)\ol g(t). \tag 8.8$$
If (8.8) is not true, there exists a sequence $t_j\to 0$ such that
$\ol g(\frac 12 t_j)\ge (\frac 12-2^{-j})\ol g(t_j)$. Let $\ol
g_j(t)=a_j^{-1}\ol g(t_jt)$ and let $u^*_j(x_1,
x_2)=a_j^{-1}u^*(a_jx_1, t_jx_2)$, where $a_j=\ol g(t_j)$. Then
$u^*_j$ sub-converges to a convex function $u^*_0$, and $\ol g_j$
converges to a linear function $\ol g_0$  and  $\ol g_0(t)=
u^*_0(0, t)+u^*_0(0, -t)$. By convexity it follows that $u^*_0(0,
t)$ is linear in both $t>0$ and $t<0$. On the other hand, since
$\lim_{x_1\to 0}\frac 1{|x_1|}u^*(x_1, 0)=1$, we have $u^*_0(x_1,
0)=|x_1|$. Hence $|N_{u^*_0}(\{0\})|>0$, that is the singular part
of the Monge-Amp\`ere measure $\mu[u^*_0]$ does not vanish at the
origin. On the other hand, since $u^*$ is a maximizer (see
(7.14)), $u^*_0$ is also a maximizer, and so the singular part of
$\mu[u^*_0]$ must vanish by the proof of Lemma 2.3. We reach a
contradiction. $\square$

Let $\alpha^\pm=\alpha^\pm_h$ be the constants such that
$$[t^+_h]^{1+\alpha^+}=h,\ \ \ \ [t^-_h]^{1+\alpha^-}=h,$$
and let $\alpha_h=\min(\alpha^+, \alpha^-)$. By (8.5) and (8.7) we
have
$$\beta\le \alpha_h\le 1.  $$

Denote $\alpha_0=\underline{\text{lim}}_{h\to 0} \alpha_h$. Choose
a sufficiently small $h>0$ such that
$$\alpha_h\le \alpha_0+\eps, \tag 8.9$$
where $\eps>0$ is a sufficiently small constant. Denote
$\alpha=\min(\alpha_h^+, \alpha^-_h)$. Suppose without loss of
generality that $\alpha =\alpha^+_h$. Denote
$\Om_h^+=\Om_h\cap\{x_2>-h^{1/(1+\alpha)}\}$ and
$q_h=|N_{u^*}(\Om_h^+)|$.  Then by Lemma 8.2,
$$|\Om^+_h|\le C h^{1+1/(1+\alpha)}. \tag 8.10$$
We need an lower bound for $q_h$.

\proclaim{Lemma 8.4} For $h>0$ small, we have
$$q_h \ge C h^{\alpha/(1+\alpha)} \tag 8.11$$
for some $C>0$ independent of $h$.
\endproclaim

\noo{\it Proof.}\ As in the proof of Lemma 8.2, we denote
$u^*_h(x_1, x_2)=h^{-1}u^*(hx_1, t_h^+x_2)$, where
$t^+_h=h^{1/(1+\alpha)}$. Then (8.11) is equivalent to
$|N_{u^*_h}( G_h)|\ge C$, where $G_h=\{x\in\R^2\ |\ u^*_h(x)<1,
x_2>-1\}$.

Let $\omega_h=\{x\in\R^2\ |\ u^*_h(x)<a(x_2+1)\}$. By (8.6) we see
that $\omega_h\subset G_h$ if $a>0$ is small (one can choose
$a=\frac 13$ for $h>0$ small). Let $z_h$ be a convex function such
that $z(0)=0$, $z(x)=u^*_h(x)$ on $\p \omega_h$, and $z(tx)=tz(x)$
for any $t\ge 0$. Then the graph of $z$ is a convex cone. It is
easy to see that
$$|N_{u^*_h}(G_h)|\ge  |N_{z_h}(\{0\})| \ge C. $$
Hence (8.11) holds. $\square$

By a rescaling we suppose $B_1(0)\subset \Om^*$. Let $v$ be the
solution of
$$\mu[v]=(1-\chi_{\Om_h^+})\mu[u^*]+\delta\chi_{_{\wtt E}}\tag 8.12$$
such that $v=u^*$ on $\p B_1$, where
$$\wtt E=\Om_{h^*}\cap\{|x_2|<\frac 12 (h^*)^{1/(1+\alpha^*)} \},$$
$h^*$ is given in (8.20) below,  and $\chi_{_{\wtt E}}$ is the
characteristic function of $\wtt E$. We choose the constant
$\delta>0$ such that
$$\delta |\wtt E|=\frac 12 q_h.\tag 8.13$$

We claim that there is a point $x_0\in \Om^+_h$ such that
$$v(x_0)>u^*(x_0)+c_0 h^{\alpha/(1+\alpha)}\tag 8.14$$
for some constant $c_0>0$ independent of $h$. Indeed, if (8.14) is
not true, let
$$\ol v=v+c_0 h^{\alpha/(1+\alpha)} (|x|^2-1).$$
Then $\ol v\le u^*$ on $\p (B_1-\Om^+_h)$ and $\mu[\ol v]\ge
\mu[v]\ge \mu[u^*]$ in $B_1-\Om^+_h$. Hence by the comparison
principle we conclude that $\ol v\le u^*$ in $B_1-\Om^+_h$. It
follows that
$$|N_{\ol v}(B_1)|\ge |N_{u^*}(B_1)|.\tag 8.15$$
The following lemma provides an upper bound for $|N_{\ol
v}(B_1)|$.

\proclaim{Lemma 8.5} We have
$$|N_{\ol v}(B_1)|\le |N_v(B_1)|
                   + c_1h^{\alpha/(1+\alpha)},\tag 8.16$$
where $c_1$ can be arbitrarily small as long as $c_0$ is
sufficiently small.
\endproclaim

To use Lemma 8.5, from (8.12) we have
$$|N_v(B_1)| = |N_{u^*}(B_1-\Om_h^+)|+\frac 12 q_h
             = |N_{u^*}(B_1)|-\frac 12 q_h .$$
Hence when $c_0>0$ is chosen small, by (8.11) and (8.16) we have
$$|N_{\ol v}(B_1)|\le |N_{u^*}(B_1)|-\frac 12 q_h
                   + c_1h^{\alpha/(1+\alpha)}<|N_{u^*}(B_1)|, $$
which contradicts (8.15). Hence (8.14) holds.

\noo{\it Proof of Lemma 8.5.}\  By approximation we may suppose
$v$ is smooth. Then $G=Dv(B_1)$ is a bounded topological disc. Let
$\p B_1$ be parametrized by $\th\in [0, 2\pi]$, namely
$x=(\cos\th, \sin\th)$. Then $Dv$ is a diffeomorphism from $\p
B_1$ to $\p G$. By convexity we have
$$\lan\frac {d}{d\th} Dv(x), \frac {d}{d\th} x\ran>0
                              \ \ \ x\in\p B_1,\tag 8.17$$
where $\lan \cdot, \cdot\ran$ denotes the inner product in $\R^2$.
For any $x\in\p B_1$, we have
$$D\ol v(x)=Dv(x)+ b x ,\tag 8.18$$
where $b=2c_0 h^{\alpha/(1+\alpha)}$. If $G$ is a ball, (8.16)
follows immediately from (8.17) and (8.18), with $c_1=4\pi c_0$.

In general we choose $R>0$ such that $G$ is contained in the ball
$B_R$. Observe that $N_{\ol v}(B_1)-N_v(B_1)$ is the region
covered by the family of line segments $\{\ell_x\}_{x\in\p B_1}$
with endpoints $Dv(x)$ and $Dv(x)+ b x$. We move these line
segments to new positions, $\ell_x\to\wtt \ell_x$, such that that
one endpoint is $z_x\in\p B_R$ and the other one is $z_x+b x$, and
for any $x$, both $\ell_x$ and $\wtt\ell_x$ are on the same
straight line. Then by (8.17), the area of the region covered by
$\{\wtt\ell_x\}$ is larger than that of the region covered by
$\{\ell_x\}$. Hence (8.16) holds with $c_1=4\pi c_0R$. $\square$

By Lemma 8.2, $\Om_h\subset \{-h\le x_1\le h\}$. Hence we may
suppose (8.14) holds at some point $x_0\in \{x_1=0\}$, namely
$x_0=(0, a_0)$. Noting that  by (8.6) and our definition of
$\alpha$,  $\Om_h^+\subset \{x_2< C h^{1/(1+\alpha)} \}$, we have
$$a_0\le C h^{1/(1+\alpha)}. $$
Denote $\eta(t)=v(0, t)$. Suppose $\eta'(a_0)\ge 0$, (the case
$\eta'(a_0)\le 0$ can be dealt with similarly). Then
$$\eta(t)>g(t)+\frac {c_0}2 h^{\alpha/(1+\alpha)} \tag 8.19$$
if $t\ge a_0$ and $g(t)<\frac 12c_0 h^{\alpha/(1+\alpha}$. Let
$$h^*=\frac 12c_0 h^{\alpha/(1+\alpha)}. \tag 8.20$$
Then (8.19) holds with
$$a_0\le t\le (h^*)^{1/(1+\alpha^*)},$$
where $\alpha^*=\alpha(h^*)$.

It follows that
$$v>u^*\ \ \ \text{in}\ \ \
         G_{h^*/2}:=\Om_{h^*/2}\cap\{x_2>a_0\}.\tag 8.21$$
Denote
$$E=\{x\in \R^2\ |\ 2(x-x_0)+x_0\in G_{h^*/2}\}, $$
namely $E$ is the $\frac 12$-dilation of $G_{h^*/2}$ with respect
to $x_0$. Then by (8.6),
$$|E|\ge C (h^*)^{1+1/(1+\alpha^*)} \ge C_1 |\wtt E| .\tag 8.22$$
We also have
$$|N_{u^*}(E)|\le C(h^*)^{\alpha^*/(1+\alpha^*)}.\tag 8.23$$
To verify (8.23) one observes that $N_{u^*}(E)$ is contained in
the set $\{x\in\R^2\ |\ |x_1|<1, |x_2|\le
C(h^*)^{\alpha^*/(1+\alpha^*)}\}$.

Let $\hat v=\max (v, u^*)$. We claim
$$A^*(\hat v)>A^*(u^*),\tag 8.24$$
which implies that $u^*$ is not a maximizer, and so we reach a
contradiction. Observing that by (8.21),  $\hat v=v$ in $E$,  and
by (8.12), $\mu[\hat v]\ge \mu[u^*]$ in $B_1(0)-(E\cup \Om_h)$, we
need only to prove
$$\int_E (\mu[v])^{3/4}>\int_E (\mu[u^*])^{3/4}
             +\int_{\Om^+_h}(\mu[u^*])^{3/4},\tag 8.25$$
where $\mu[v]=\det \p^2 v$ and $\mu[u^*]=\det \p^2 u^*$, see
(7.3). By Lemma 2.1 and 2.3, $\mu[u^*]$ is an integrable function.

Let $E=E_1\cup E_2\cup \Om_h^+$, where $E_1=\{x\in E-\Om^+_h\ |\
\mu[u^*](x)<\frac 12 \delta\}$ and $E_2=E-(E_1\cup\Om_h^+)$. If
$|E_1|\ge \frac 12 |E|$, then for any $x\in E_1$, we have
$$(\mu[v])^{3/4} \ge (\mu[u^*])^{3/4}+C\delta^{3/4}. $$
On the other hand,
$$ \int_{\Om^+_h} (\mu[u^*])^{3/4}
          \le |\Om_h^+|^{1/4}
          \big(\int_{\Om^+_h}\mu[u^*]\big)^{3/4}
          =  q_h^{3/4}|\Om_h^+|^{1/4}. \tag 8.26$$
Hence we need only
$$\delta^{3/4}|E_1|\ge Cq_h^{3/4}|\Om_h^+|^{1/4},$$
namely $|E|\ge C|\Om_h^+|$ by (8.13) and (8.22), which is
obviously true.

Otherwise observing that $|\Om_h^+|=o(|E|)$ as $h\to 0$, we may
suppose simply that $E_1$ is an empty set by replacing $E$ by
$E_2$. Then at any point $x\in E$,
$$\align
\mu[v]^{3/4}
 & =(\mu[u^*]+\delta)^{3/4}\\
 & \ge (\mu[u^*])^{3/4}+C(\mu[u^*])^{-1/4}\delta.\\
 \endalign $$
By the H\"older inequality we have
$$\int_E (\mu[u^*])^{-1/4}\ge |E|^{5/4}
       \big[\int_E \mu[u^*]\big]^{-1/4}.$$
Hence
$$\align
\int_E (\mu[v])^{3/4}
  & \ge \int_E (\mu[u^*])^{3/4}
  +C\delta |E|^{5/4}\big[\int_E \mu[u^*]\big]^{-1/4}\\
  & = \int_E (\mu[u^*])^{3/4}+C q_h |E|^{1/4}
          \big[\int_E \mu[u^*]\big]^{-1/4}.\\
          \endalign $$
Therefore by (8.26) we need only
$$q_h^{1/4} |E|^{1/4}>
    C |\Om_h^+|^{1/4} \big[\int_E \mu[u^*]\big]^{1/4} .$$
By (8.10), (8.11), (8.22) and (8.23), we need only
$$h^{\frac {\alpha}{1+\alpha}}(h^*)^{1+\frac 1{1+\alpha^*}}
  \ge C h^{1+\frac 1{1+\alpha}}(h^*)^{\frac
  {\alpha^*}{1+\alpha^*}}.$$
By (8.20), the above inequality holds if $\alpha<1+\alpha^*$,
which is true by our choice of $h$ in (8.9).  This completes the
proof of Theorem 8.2. $\square$

\vskip10pt

Next we extend Theorem 8.2 to locally convex surfaces.

\proclaim{Theorem 8.3} Suppose $n=2$. Then the maximizer in
Theorem 3.1 is strictly convex.
\endproclaim

\noo{\it Proof}. Let $\M$ denote the maximizer. If $\M$ contains a
line segment $\ell$, let $P$ be a tangent plane of $\M$ in $\R^3$
which contains $\ell$, and $F$ the component of $P\cap\M$ which
contains $\ell$. Then $F$ is a convex set [22]. If $F$ has an
extreme point which is an interior point of $\M$, then Case a) in
the proof of Theorem 8.2 applies and we reach a contradiction.

If all extreme points of $F$ are boundary points of $\M$, we are
in Case b) of the proof of Theorem 8.2. In this case, by a
rotation of axes we may suppose $F\subset\{x_3=0\}$. Then the
south pole of the sphere $S^2$ is an interior point of the Gauss
mapping image of $\M_0$. Let $F_\delta$ denote the connected
component of $\M\cap\{x_3<\delta\}$ containing $F$. Then
$F_\delta$ is convex if $\delta>0$ is sufficiently small [22].

Let $F_\delta^-$ denote the set of points $p\in F_\delta$ such
that the Gauss mapping $G(p)$ of $\M$ falls in the south
hemisphere. Then $F^-_\delta$ is the graph of an affine maximal
convex function $u$. Since all extreme points of $F$ are boundary
points of $\M$, there is a line segment $\ell\subset F$, as given
in (8.1), of which both endpoints $\pm e_1$ are boundary points of
$\M$. We claim that $\ell$ is transversal to $\p\M$ at the
endpoints $\pm e_1$. Indeed, extend $\M_0$ to $\wtt\M_0$ such that
$\wtt\M_0$ is locally uniformly convex and
$\M_0\cap\p\wtt\M_0=\emptyset$. Denote $\wtt\M=\M\cup(\wtt
M_0-\M_0)$. Then $\wtt\M$ is a locally convex hypersurface. For
$\delta>0$ we denote by $\wtt F_\delta$ the component of $\wtt
\M\cap\{x_3<\delta\}$ containing $\ell$. Then $\wtt F_\delta$
contains no boundary points of $\wtt\M$ when $\delta>0$ is small.
Hence by Lemma 3.1, $\wtt F_\delta$ is convex and so $\ell$ is
transversal to $\p\M$.

It follows that $u$ satisfies (8.2). Let $u^*$ be the Legendre
transform of $u$. Then $u^*$ satisfies (8.3) and (8.4). Hence the
proof of Theorem 8.2 applies and we also reach a contradiction.
$\square$

To conclude this section we extend Theorem 8.2 to the
inhomogeneous equation (4.1) for general $f\in L^\infty(\Om)$. We
will need the following lemma.

\proclaim{Lemma 8.6} Let $\Om$ be a uniformly convex domain and
$\phi_k$ be a sequence of convex functions converging uniformly to
$\phi$. Suppose $u_k\in \ol S[\phi_k, \Om]$ converges to $u_0$ and
$$A_\alpha (u_k)\ge
    \sup_{v\in \ol S[\phi_k, \Om]} A_\alpha (v)-\eps_k \tag 8.27$$
with $\eps_k\to 0$ as $k\to\infty$. Then
$$A_\alpha (u_0)=\sup_{v\in \ol S[\phi, \Om]} A_\alpha (v),\tag 8.28$$
where $A_\alpha(u)=\int_\Om [\det\, \p^2 u]^\alpha$ and $\alpha\in
(0, 1)$ is a constant.
\endproclaim

\noo{\it Proof.}\ If $\phi_k=\phi$ for all $k$, Lemma 8.6 follows
from the upper semi-continuity of the functional $A_\alpha$.

In the general case we choose two sufficiently small constants
$\delta>>\eps>0$ and a function $\ol\phi\in C^2(\pom)$ with $0\le
\ol\phi-\phi_k\le \eps$ for sufficiently large $k$, such that for
any function $u\in\ol S[\phi_k, \Om]$, we can modify $u$ in the
open set $\{x\in\Om;\ \dist(x, \pom)<\delta\}$ such that
$u=\ol\phi$ on $\pom$.  Then Lemma 8.6 follows from the
equi-integrability of $[\det \p^2 u]^\alpha$ for $\alpha\in
(0,1)$, see the proof of Lemma 7.1. $\square$

\proclaim{Theorem 8.4} Let $\Om$ and $\phi$ be as Theorem 8.2.
Suppose $n=2$ and $f\in L^\infty(\Om)$. Then the maximizer $u$ in
Theorem 2.1 is strictly convex.
\endproclaim

To prove Theorem 8.4 we will use rescaling to reduce to the
homogeneous case $f\equiv 0$ and apply the proof of Theorem 8.2.
Let $y=T(x)$ be a linear transformation and $v(y)=Mu(x)$ for some
constant $M>0$, so that $v$ satisfies
$$L(v)= \beta L(u) , \tag 8.29$$
where $\beta=(M |T|)^{-2/(n+2)}$. If $u$ is a maximizer of $A_f$
over the set $\ol S[\phi, \Om]$, where $A_f$ is given in (2.18),
then $v$ is a maximizer of the functional
$$A_{\beta f}(u)=\int_{\hat \Om} [\det \p^2 u]^{1/(n+2)}
          -\beta \int_{\hat \Om} fu \tag 8.30$$
over the set $\ol S[\hat \phi, \hat \Om]$, where $\hat
\phi(y)=M\phi(x)$ and $\hat \Om=T(\Om)$.

The linear transformation $T$ will be chosen to normalize a given
bounded convex domain. For a given bounded convex domain
$\Om\subset \R^n$, there is a unique ellipsoid $E$, called the
minimum ellipsoid of $\Om$, which has the minimum volume among all
ellipsoids containing $\Om$, such that
$$\frac 1n E\subset\Om\subset E, \tag 8.31$$
where $\alpha E$ is the concentrated $\alpha$-dilation of $E$. We
say $\Om$ is normalized if $E$ is the unit ball. For a general
convex domain $\Om$ we may choose a linear transformation $T$ such
that $T(E)$ is the unit ball and $T(\Om)$ is normalized.

\vskip5pt

\noo{\it Proof of Theorem 8.4}.\ If $u$ is not strictly convex,
then the graph of $u$, $\M_u$, contains a line segment $\ell$. Let
$P$ denote the tangent plane of $\M_u$ containing $\ell$, and $F$
the component of $P\cap\M_u$ which contains $\ell$. By adding a
linear function to $u$, we may suppose $P=\{x_3=0\}$. Then $F$ is
a convex set in the plane $\{x_3=0\}$.

If $F$ contains an extreme point which is an interior point of
$\M_u$, by a rotation of the coordinates we suppose the extreme
point is the origin and $F\subset \{x_1\ge 0\}$ such that
$\{x_1=0\}\cap\ol F$ contains the origin only. Without loss of
generality we assume that $\ol F\cap\{x_1<1\}$ is strictly
contained in $\Om$.  By a linear transformation we also suppose
the line segment $\{(t, 0)\ {|}\ \ 0\le t\le 1\}$ is contained in
$F$. Let $\ell_\eps (x) = \eps (1- x_1)$ be a linear function such
that $\ell_\eps>u$ near the origin and the domain
$\omega_\eps=:\{\ell_\eps>u\}$ is strictly contained in $\Om$. Let
$v_\eps (y)=\eps^{-1}(\ell_\eps-u)(x)$, where $y=T_\eps(x)$,
$x\in\omega_\eps$, and $T_\eps$ is a linear transformation
normalizing the domain $\omega_\eps$. Then $v_\eps\le 0$ in
$T_\eps(\omega_\eps)$ and by (8.29), we have $L(v_\eps) =
\delta_\eps f$ with $\delta_\eps = (\eps/|T|)^{2/(n+2)}\to 0$ as
$\eps\to 0$. Let $\M_\eps$ denote the graph of $v_\eps$. Then
$\M_\eps\subset \{x_{n+1}\le 0\}$ and it converges to a locally
convex hypersurface $\M_0$ as $\eps\to 0$. By Lemma 8.6, $\M_0$ is
affine maximal. Note that to use Lemma 8.6, one may need to
decompose $\M_0$ into finitely many pieces such that each piece is
a graph with uniformly bounded gradient in appropriate coordinate
systems. Since the boundary of $\M_0$ lies on the plane
$\{x_3=0\}$, we conclude that $\M_0$ is strictly convex, as noted
after Theorem 8.1. But on the other hand, $\M_0$ contains a line
segment by our construction, a contradiction. Hence no extreme
point of $F$ can be interior point of $\M_u$.

It follows that all extreme points of $F$ are boundary points of
$\M_u$. In this case there exists a line segment, which we suppose
is $\ell=\{te_1\ \big{|}\ -1\le t\le 1\}$, such that both
endpoints of $\ell$ are boundary points of $\M_u$. Moreover the
origin is an interior point of $\Om^*=D\phi(\Om)$, namely there is
$R>0$ such that the disc $B_R(0)$ is contained in $\Om^*$. Let
$u^*$ be the Legendre transform of $u$. Then $u^*(0)=0$ and (see
Remark 8.4)
$$u^*(x)\ge x\cdot x_0\ \ \ \forall\ \ x_0\in F.\tag 8.32$$

For $h>0$ small, denote $\Om_h=\{x\in \Om\ {|}\ u(x)<h\}$ and
$u_h(y)=h^{-1}u(x)$, where $y=T_h(x)$ and $T_h(x_1, x_2)=(x_1,
x_2/b_h)$ is a linear transformation with
$$b_h=\sup\{|x_2|:\ (x_1, x_2)\in \Om_h\}.$$
Then the set $T_h(\Om_h)$ is uniformly bounded and the area
$|T_h(\Om_h)|$ has a positive lower bound. Let $\M_h$ denote the
graph of $u_h$. By choosing a subsequence, we may suppose $\M_h$
converges to $\M_0$ as $h\to 0$. By Lemma 8.6, $\M_0$ is affine
maximal. Moreover, the line segment $\ell$ is contained in $\M_0$
with both its endpoints lying on $\p \M_0$.

Let $u^*_h$ denote the Legendre transform of $u_h$. We have
$$u^*_h(x_1, x_2)=\frac 1h u^*(hx_1, hx_2/b_h).$$
The function $u^*_h$ is uniformly bounded in $B_\rho (0)$ for some
$\rho>0$, with an upper bound independent of $h$. This is because
for any plane $x_3=a\cdot x-b$ of slope $\rho$, i.e. $|a|=\rho$,
with $b$ sufficiently large, the plane lies below and does not
touch the graph of $u_h$. Let $b_0$ be the least constant with
this property. Then $u^*_h(a)=b_0$.  Hence $\sup_{\p B_\rho(0)}
u^*_h$ is bounded. By convexity, $u^*_h$ is uniformly Lipschitz in
$B_{\rho/2}(0)$.

It follows that for any sequence $h_k\to 0$ (which will be chosen
such that (8.34) holds), there is a subsequence such that
$u_{h_k}^*$ converges locally to a convex function $u_0^*$. From
Lemma 8.6 (with $\Om=B_\rho(0)$ and $\phi_k=u_{h_k}$) and Lemma
7.1, $u_0^*$ is a maximizer of the functional $A^*$.

Since the line segment $\ell$ lies in the graph of $u_h$ for any
$h>0$, we have $u_h^*(x)\ge |x_1|$ for any $h>0$. It follows
$u^*_0(x)\ge |x_1|$ in $B_\rho(0)$. From Remark 8.4, $\lim_{t\to
0} \frac 1t u^*(t, 0)=1$. Hence we have
$$u^*_0(t, 0)=|t|. \tag 8.33$$

To use the proof of Theorem 8.2, we need to show that if the
subsequence $\{h_k\}$ is properly chosen, then the function $g_0$,
where $g_0(t)=u^*_0(0, t)$, is of polynomial growth near $t=0$.
Namely there exist positive constants $\beta\ge \alpha>0$ such
that
$$|t|^{1+\beta}\le g_0(t)\le |t|^{1+\alpha}\tag 8.34$$
near $t=0$. Once (8.34) is proved, the proof of Theorem 8.2
applies and so $u$ must be strictly convex.

Let $g(t)=u^*(0, t)$. Then $g$ satisfies (8.5). From the proof of
Lemma 8.3, we see that $g$ satisfies (8.8). Hence $g_0$ also
satisfies (8.8), namely $\ol g_0(\frac 12 t)\le (\frac 12-\th)\ol
g_0(t)$ for any $t>0$, where $\ol g_0(t)=g_0(t)+g_0(-t)$. Hence
the second inequality in (8.34) holds.

The first inequality of (8.34) follows from (8.5). Indeed, we
claim that for any $t_0>0$, there is a $t_1<t_0$ such that for any
$t\in (0, \frac 12t_1)$,
$$g(t)\ge (t/t_1)^3g(t_1).\tag 8.35$$
If the claim is true, then there is a sequence $t_k\searrow 0$
such that $g(t)\ge (t/t_k)^3g(t_k)$ for any $t\in (0, \frac 12
t_k)$. Hence the limit function $g_0$ satisfies $g_0(t)/g_0(1)\ge
t^3$ $\forall$ $t\in (0, \frac 12)$. Namely first inequality in
(8.34) holds.

To prove the claim we suppose on the contrary that (8.35) is not
true. Then there exists $t_1=\th_1 t_0$ with $\th_1\in (0, \frac
12)$ such that $g(t_1)\le \th_1^3 g(t_0)$. Also there exists
$t_2=\th_2 t_1$ with $\th_2\in (0, \frac 12)$ such that $g(t_2)\le
\th_2^3 g(t_1)$. Repeat this process we obtain a sequence
$\{t_k\}$,  $t_k=\th_kt_{k-1}$ with $\th_k\in (0, \frac 12)$, such
that
$$\align
g(t_k) & \le g(t_{k-1})\, \th_k^3\le \cdots\\
    & \le g(t_0) [\th_1\cdots\th_k]^3
         = \big[g(t_0)\frac {t_k}{t_0^3}\big]\, t_k^2.\\
     \endalign $$
When $k$ is sufficiently large, $g(t_0)\frac {t_k}{t_0^3}<1$, we
reach a contradiction with (8.5). Hence (8.34) holds. This
completes the proof.  $\square$

Combining Theorem 8.4 with Theorem 6.2 and Theorems 4.1 and 4.2,
we therefore obtain

\proclaim{Theorem 8.5} Let $\Om$ be a bounded smooth domain in
$\R^2$, and $\phi\in C^2(\bom)$ be a uniformly convex function.
Suppose $f\in L^\infty(\Om)$. Then there exists a unique, locally
uniformly convex maximizer $u\in W^{4, p}_{loc}(\Om)$ for (2.17).
If furthermore $f\in C^{k, \alpha}(\Om)\cap L^\infty(\Om)$,
$\alpha\in (0, 1)$, then $u\in C^{k+4, \alpha}(\Om)$.
\endproclaim

\vskip10pt

\baselineskip=14.5pt
\parskip=1pt

\centerline{\bf \S 9. Final remarks}

\vskip10pt

\noo{\it Remark 9.1}. \ Let $u$ be a smooth maximizer in Theorem
B. Obviously $u$ satisfies the first boundary condition (1.9) and
$Du(\Om)\subset D\phi(\bom)$. An interesting problem is whether
$u$ satisfies the second boundary condition $Du(\Om)=D\phi(\Om)$.
Recall that the Dirichlet problem of the minimal surface equation
is solvable for any smooth boundary values if and only if the
boundary is mean convex. Analogously  for the boundary condition
(1.9), additional conditions may be necessary in order that
$Du(\Om)=D\phi(\Om)$ holds.

In [19] we studied the regularity of hypersurfaces $\M_\Ga$
contained in a given closed, uniformly convex hypersurface $\Ga$
with maximal affine area, and proved that the contact set
$\M_\Ga\cap\Ga$ must be affine mean convex. Therefore a reasonable
assumption such that $Du(\Om)=D\phi(\Om)$ holds is that $\phi$ is
affine mean convex, namely the affine mean curvature of the graph
of $\phi$, $H_A[\phi] = -\frac {1}{n+1} L[\phi]$, is positive.

\noo{\it Remark 9.2}.  In [21] we proved that the function $u$ in
$\R^{10}$, given by
$$u(x)=\sqrt {|x'|^9 + x_{10}^2} , $$
where $x'=(x_1, \cdots, x_9)$, is affine maximal. Therefore in
high dimensions $(n\ge 10)$ there exist affine maximal functions
vanishing on the boundary of a convex domain which are not
strictly convex. We have not found similar examples in dimensions
$3\le n\le 9$.

\noo{\it Remark 9.3}. Let $u$ be a convex solution to the
Monge-Amp\`ere equation
$$\det D^2 u=1.$$
Then $u$ satisfies the affine maximal surface equation (1.5). When
$n\ge 3$, there exist non-smooth convex solutions to the above
Monge-Amp\`ere equation [17]. Hence there exist non-smooth affine
maximal hypersurfaces when $n\ge 3$.

However a convex solution to the general Monge-Amp\`ere equation
$$\det D^2 u=f\ \ \ \text{in}\ \Om,$$
where $f$ is a bounded, positive function, is strictly convex in
any dimension if $u$ vanishes on $\pom$ or both $\pom$ and the
trace of $u$ on $\pom$ are  smooth. For the affine maximal surface
equation, an interesting question is whether the maximizer $u$ in
Theorem $B$ is strictly convex if $\phi$ is uniformly convex.

 \vskip30pt

\baselineskip=12pt
\parskip=1pt

\Refs\widestnumber\key{ABC}

\vskip10pt

\item {[1]}   W. Blaschke,
              Vorlesungen \"uber Differential geometrie, Berlin, 1923.

\item {[2]} L.A. Caffarelli,
             A localization property of viscosity solutions to
             the Monge-Amp\`ere equation and their strict convexity,
             Ann. Math., 131(1990), 129-134.

\item {[3]} L.A. Caffarelli,
             Interior $W^{2,p}$ estimates for solutions
             of Monge-Amp\`ere equations,
             Ann. Math., 131(1990), 135-150.

\item {[4]} L.A. Caffarelli and C.E. Guti\'errez,
             Properties of the solutions of the linearized
             Monge-Amp\`ere equations,
             Amer. J. Math., 119(1997), 423-465.

\item {[5]} L.A. Caffarelli, L. Nirenberg, and J. Spruck,
             The Dirichlet problem for nonlinear second order
             elliptic equations I. Monge-Amp\`ere equation,
             Comm. Pure Appl. Math., 37(1984), 369-402.

\item {[6]} E. Calabi,
              Hypersurfaces with maximal affinely invariant area,
              Amer. J. Math. 104(1982), 91-126.

\item {[7]} E. Calabi,
              Affine differential geometry and holomorphic curves,
              Lecture Notes Math. 1422(1990), 15-21.

\item {[8]} S.Y. Cheng and S.T. Yau,
             Complete affine hypersurfaces, I.
             The completeness of affine metrics,
             Comm. Pure Appl. Math., 39(1986), 839-866.

\item {[9]} S.S. Chern,
             Affine minimal hypersurfaces,
             in {\it minimal submanifolds and geodesics},
             (Proc. Japan-United States Sem., Tokyo, 1977, 17-30.

\item {[10]} D. Gilbarg and N.S. Trudinger,
             Elliptic partial differential equations of second order,
             Springer-Verlag, New York, 1983.

\item {[11]} B. Guan and J. Spruck,
             Boundary value problems on $S\sp n$ for surfaces of
             constant Gauss curvature. Ann. of Math. 138 (1993),
             601--624.

\item {[12]}  N.Ivochkina,
             A priori estimate of $\|u\|_{C^2(\bom)}$ of convex
             solutions of the Dirichlet problem for the Monge-Amp\`ere
             equation.
             Zap. Nauchn. Sem. Leningrad. Otdel. Mat. Inst. Steklov.
             (LOMI) 96 (1980), 69--79 (Russian). English
             translation in J. Soviet Math., 21(1983), 689-697.

\item {[13]}   N.V. Krylov,
              Nonlinear elliptic and parabolic equations of the
              second order, Reidel, Dordrecht-Boston, 1987.

\item {[14]}   K. Leichtweiss,
              Affine geometry of convex bodies,
              Johann Ambrosius Barth Verlag, Heidelberg, 1998.

\item {[15]}   E. Lutwak,
              Extended affine surface area,
              Adv. Math., 85(1991), 39-68.

\item {[16]}  K. Nomizu and T. Sasaki,
              Affine differential geometry,
              Cambridge University Press, 1994.

\item {[17]}    A.V. Pogorelov,
               The muitidimensional Minkowski problems,
               J. Wiley, New York, 1978.

\item {[18]}   M.V. Safonov,
              Classical solution of second-order nonlinear elliptic
              equations, Izv. Akad. Nauk SSSR Ser. Mat. 52 (1988),
              1272--1287 (Russian).  English translation in
              Math. USSR-Izv. 33 (1989), 597--612.

\item {[19]}  W.M. Sheng, N.S. Trudinger, and X.-J. Wang,
             Enclosed convex hypersurfaces with maximal affine
             area, preprint.

\item {[20]}   U. Simon,
              Affine differential geometry,
              in {\it Handbook of differential geometry},
              North-Holland, Amsterdam, 2000, 905-961.

\item {[21]} N.S. Trudinger and X.-J. Wang,
              The Bernstein problem for affine maximal hypersurfaces,
              Invent. Math., 140 (2000), 399--422.

\item {[22]} N.S. Trudinger and X.-J. Wang,
              On locally convex hypersurfaces with boundary,
              J. Reine Angew. Math., 551(2002), 11-32.

\item {[23]} N.S. Trudinger and X.-J. Wang,
             On boundary regularity for Monge-Amp\`ere equation
             and applications, in preparation.

\item {[24]}   X.-J. Wang,
             Affine maximal hypersurfaces,
             Proc. ICM  Vol.I\!I\!I, 2002, 221-231.

\item {[25]}   W. P. Ziemer,
              Weakly differentiable functions,
              Springer, 1989.

\endRefs

\enddocument

\end